\documentclass[12pt,a4paper,english]{article}
\usepackage[affil-it]{authblk}
\usepackage[utf8]{inputenc}
\usepackage[margin=1in]{geometry}
\usepackage[english]{babel}
\fontfamily{phv}\selectfont
\usepackage{stfloats}
\usepackage{cite}
\usepackage{amsthm,amsmath,amssymb,amsfonts,mathtools}
\usepackage{longtable}
\usepackage{subfigure}
\usepackage{graphicx}
\usepackage{textcomp}
\usepackage[hyphens]{url}
\usepackage[hidelinks,breaklinks]{hyperref}

\theoremstyle{remark}

\theoremstyle{definition}
\newtheorem{definition}{Definition}

\begin{document}
\title{\textbf{Optimal control of a commercial building's thermostatic load for off-peak \\demand response}}
\author{Randall Martyr\thanks{Corresponding author. Email: r.martyr@qmul.ac.uk} \thanks{Financial support received from the UK Engineering and Physical Sciences Research Council (EPSRC) via Grant EP/N013492/1.}}
\author{John Moriarty \thanks{Financial support received from the UK Engineering and Physical Sciences Research Council (EPSRC) via Grant EP/P002625/1.}}
\author{Christian Beck\protect\footnotemark[2]}
\affil{School of Mathematical Sciences, Queen Mary University of London, Mile End Road, London E1 4NS, United Kingdom.}
\maketitle
\begin{abstract}
	This paper studies the optimal control of a commercial building's thermostatic load during off-peak hours as an ancillary service to the transmission system operator of a power grid. It provides an algorithmic framework which commercial buildings can implement to cost-effectively increase their electricity demand at night while they are unoccupied, instead of using standard inflexible setpoint control. Consequently, there is minimal or no impact on user comfort, while the building manager gains an additional income stream from providing the ancillary service, and can benefit further by pre-conditioning the building for later periods. The framework helps determine the amount of flexibility that should be offered for the service, and cost optimized profiles for electricity usage when delivering the service. Numerical results show that there can be an economic incentive to participate even if the payment rate for the ancillary service is less than the price of electricity.
\end{abstract}


\section*{Notation}

\begin{longtable}{cp{.60\textwidth}c}
	\centering
	& \centering {\bf Constant quantities} & \\
	{\bf Notation} & {\bf Description} & {\bf Units} \\
	$T$ & length of control horizon, the period of time during which the building's temperature is controlled for the ancillary service & min \\
	$P$ & night-time electricity price & p/kWh \\
	$R$ & utilization payment & p/kWh \\
	$X_{min}$ & night-time lower temperature limit & $\,^{\circ}\rm{C}$ \\
	$X_{max}$ & night-time upper temperature limit & $\,^{\circ}\rm{C}$ \\
	$\hat{X}$ & temperature limit used for pre-cooling at time $T$ & $\,^{\circ}\rm{C}$ \\
	$X_{off}$ & asymptotic temperature for the building ``off'' state & $\,^{\circ}\rm{C}$ \\
	$X_{on}$ & asymptotic temperature for the building ``on'' state & $\,^{\circ}\rm{C}$ \\
	$\tau$ & thermal time constant & min \\
	$C_{max}$ & maximum power limit for the building & kW \\
	$\mathcal{U}$ & set of normalized cooling power usage variables $u$ & -- \\[5pt]
	& \centering {\bf Time-varying quantities} & \\
	{\bf Notation} & {\bf Description} & {\bf Units} \\
	$x(t)$ & building internal temperature at time $t$ & $\,^{\circ}\rm{C}$ \\
	$C(t)$ & power usage of the cooling equipment & kW \\
	$C_{ref}(t)$ & reference cooling power usage & kW \\
	$C_{alt}(t)$ & alternative power usage used to calculate the level of reserve & kW \\
	$C_{cap}(t)$ & level of reserve capacity & kW \\
	$C_{ask}(t)$ & reserve service instructions, a profile of additional power usage that must be delivered & kW \\
	$C_{del}(t)$ & power usage when delivering the reserve service according to the instructions $C_{ask}$ & kW \\
	$u(t)$ & normalized power cooling usage & $1$ \\
	$u_{ref}(t)$ & normalized reference cooling power usage & $1$ \\
	$u_{alt}(t)$ & normalized alternative cooling power usage & $1$ \\
	$u_{del}(t)$ & normalized delivery cooling power usage & $1$ \\
\end{longtable}

\section{Introduction}
\subsection{The need for electricity balancing and ancillary services}
Electricity supply and demand on a power system must be balanced continuously in real time to ensure its stability. The system operator, an independent entity that manages the transmission system \cite[p.~3]{Kirschen2004}, balances the power system by:
\begin{itemize}
	\item increasing generation or reducing demand when there is a shortfall in supply;
	\item decreasing generation or increasing demand when there is surplus power.
\end{itemize}
The latter, which we refer to as {\it decremental} actions \cite{Szabo2017}, are increasingly relevant for power systems with high levels of intermittent generation from renewable energy sources \cite{Rothleder2014}. In order to carry out its balancing duties, the system operator procures a variety of {\it ancillary services} from third-party companies \cite[p.~106]{Kirschen2004}. There is much interest in enabling electricity consumers to provide ancillary services \cite{Hirst1998, Vardakas2015, Paterakis2017}, particularly commercial buildings due to the large flexible demand from their heating, ventilation and air conditioning (HVAC) systems \cite{Olivieri2014, Pavlak2014, Lawrence2016, Kim2016, DeConinck2016, Blum2017}. {\it Replacement reserves}, which are given more time to respond and are used as back up for faster acting services, can be most suitable in this case \cite[p.~32]{Hirst1998}.

This paper is a quantitative study of the potential for a commercial building with flexible thermostatic load to participate in a decremental replacement reserve (DRR) initiative that is modelled after a real-world example. The setting of this paper is novel in comparison to previous papers such as \cite{DeConinck2016,Blum2017}. In the present work the reserve provider bids a schedule of their reserve capacity (in kW) together with a fixed utilization price (per kWh), rather than a variable price depending on the quantity utilized. We provide a mathematical and computational framework that maximizes the benefit gained from participating in the DRR service. We also focus on temperature cooling only and note that heating can be treated symmetrically. Data centres, in particular, are an important example since they account for more than 1\% of global electricity usage, and their cooling infrastructure typically accounts for about 40\% to 50\% of their electricity usage \cite{Dayarathna2016}.

\subsection{Buildings as ancillary service providers}
{\it Demand response} refers to any programme that motivates changes in an electricity consumer's normal power usage, typically in response to incentives regarding electricity prices \cite{Vardakas2015}. It is widely considered as a cost-effective and reliable solution for improving the efficiency, reliability, and safety of the power system \cite{Vardakas2015,Paterakis2017}. There are several examples of initiatives that incentivize electricity consumers to reduce their demand, especially during peak hours, and \cite{Olivieri2014} recently studied the potential for buildings to use their HVAC systems to participate in such initiatives. Demand response schemes that provide incentives for increased electricity demand are much rarer. The Demand Turn Up (DTU) programme offered by National Grid UK, the transmission system operator in Great Britain, is one such initiative that is meant to incentivize large electricity consumers to increase their demand when there is low overall demand on the network and high output from renewable generation \cite{DemandTurnUp}. It is particularly relevant during the off-peak, night-time hours of interest to this paper, and below we summarize its key aspects (see \cite{DemandTurnUp} for further details).
\subsubsection{Demand turn up: an off-peak demand response scheme}
DTU runs during the British Summer Time (BST) period and there are two routes to market for candidates:
\begin{itemize}
	\item Fixed DTU is a medium to long-term procurement process that takes place months in advance of BST.
	\item Flexible DTU is a rolling short-term procurement process that takes place during BST and closer to the period that requires the service.
\end{itemize}
Flexible DTU can be advantageous since it gives candidates the flexibility to adjust their declaration in response to weather and market conditions. Our study is more relevant to flexible DTU since we use a dynamic model of temperature evolution that is more appropriate over the short-term.

DTU candidates declare their availability by specifying a schedule for the adjustment in electricity usage or generation they can provide, including the payment for utilization of their service. In addition to the utilization payments, successful candidates receive guaranteed payments for being available during certain windows. Unsuccessful candidates do not receive these guaranteed payments, but have the option to participate in DTU for utilization payments only.

National Grid UK sends a contracted DTU provider instructions for the service according to the capability that was declared. The provider has a deadline for acknowledging receipt of the DTU instruction, then a delivery period for responding as instructed. A DTU provider which has declared its availability must be able to deliver the service as instructed or face a penalty.

Settlement is the process of compensation for successful provision of the DTU service. A provider has two options for settlement, forecast or baseline, and its choice is fixed for the contract's duration. Both options produce a {\it reference profile} to which the actual metered electricity usage or generation is compared, and the difference is settled against the DTU service instruction that was sent. The forecast method uses the provider's prediction for electricity usage or generation during that service period, whereas the baseline method uses the average metered output from previous entries for that period in which the provider did not render a DTU service.

\subsubsection{Quantifying a building's potential to provide an ancillary service}
As mentioned previously, commercial buildings have significant potential to provide ancillary services to the transmission system operator. According to \cite[p.~1266]{Blum2017}, {\it capacity} and {\it performance} are the two main components of this ancillary service provision, each of which has a magnitude and cost. Capacity refers to the capability that the building's HVAC system has to provide an ancillary service, whereas performance refers to the work that the HVAC system does to provide the ancillary service in response to the system operator's instructions. Several papers have argued that buildings can be incentivized to participate in ancillary services markets, despite high energy prices or less efficient operating conditions, provided they are adequately compensated \cite{Olivieri2014, Pavlak2014, Lawrence2016, DeConinck2016, Blum2017}.

Reference \cite{Blum2017} proposes a methodology for quantifying flexibility and opportunity costs arising from the provision of ancillary services by buildings' HVAC systems. The authors identify sources of these opportunity costs, and develop a method of accounting for them through time that is consistent with current practice for generators. This is done by recognizing the impacts that intra-hour consumption modification associated with ancillary service provision have on daily energy efficiency and costs. The authors of \cite{DeConinck2016} previously addressed a similar problem, but unlike \cite{Blum2017} they focused on a building's capability to alter its total energy use over a period of time, and the methodology put forward was not intended for real-time dynamic operations \cite[p.~654]{DeConinck2016}. Both papers use optimal control to determine the building's capability to provide a given level of reserve and the associated opportunity cost. By varying the level of reserve, an opportunity cost curve can be constructed and used in the ancillary services market for the purpose of dispatch by the system operator, or for bidding purposes by the building manager.

\subsection{Aim of this work}
In this paper we study the potential for a commercial building to participate in an ancillary service scheme such as Demand Turn Up by controlling the electricity it consumes for temperature cooling during the night. Unlike the setting studied in \cite{DeConinck2016,Blum2017}, the reserve provider bids a schedule of their reserve capacity (in kW) together with a fixed utilization price (per kWh), rather than a variable price depending on the quantity utilized. Our main contribution is an analysis of the building manager's incentives in this novel setting by using the {\it benefit-cost ratio} of utilization payment (benefit) to night-time price of electricity (cost). By varying this ratio, we can see how it affects the magnitude of capacity offered for the ancillary service.

We approach the overall problem of economically providing the ancillary service by breaking it up into three smaller problems:
\begin{enumerate}
	\item {\it Reference:} determine an optimal reference power profile to use for settlement of the service.
	\item {\it Capacity:} determine an optimal capacity power profile to use for declaration of availability for the service.
	\item {\it Delivery:} determine an optimal delivery power profile that fulfils the service instructions.
\end{enumerate}
We formulate each of these problems as a constrained optimal control problem \cite{Clarke2013}, and use the control parametrization method \cite{Teo1991,Rehbockt1999} to obtain approximate numerical solutions for different scenarios. Optimal control is one of several mathematical techniques that can be used to optimize the provision of an ancillary service from a commercial building \cite{Deng2015, Wang2008, Olivieri2014}. Moreover, it can be an effective solution for control of the building's thermostatic load \cite[p.~158]{Vardakas2015}.

Our methodology is suitable for assessing the building's ability to participate in any demand response scheme, for either incremental or decremental reserve, where the utilization payment is fixed and the reserve provider bids capacity curves. We are able to identify the incentives that drive the optimal actions, leading to recommendations that are intuitive and implementable using a variety of control architectures. Consistent with previous studies \cite{Olivieri2014, Pavlak2014, Lawrence2016, DeConinck2016, Blum2017}, we find that, besides the dynamics and constraints for the internal temperature, the level of participation in the ancillary service depends on how well the building manager is compensated relative to the additional cost incurred. Moreover and counter-intuitively, we show that there is an economic incentive to participate even when the utilization payment is less than the night-time price of electricity.

In the following section we present our mathematical framework for optimizing the reference, capacity and delivery power profiles for the reserve service. This framework uses the building's internal temperature as a controlled variable. In principle, any model that describes the temperature dynamics using ordinary differential equations can be used. For realistic applications, the temperature relaxation behaviour of a given building is measured and used as input to the optimal control scheme. In Section~\ref{Section:Ensemble-Averaging} we consider a linear model for the temperature dynamics, which is used to obtain the numerical solutions to the optimization problems presented in Section~\ref{Section:Numerical-Results}. The paper concludes with a summary of the main results and practical recommendations in Section \ref{Section:Conclusion}.

\section{Optimal control problems for off-peak demand response from thermostatic load}\label{Section:Optimal-Control-Statement}
The {\it control horizon} is a period of time during which the building's internal temperature is controlled for the ancillary service. Let $T > 0$ denote the control horizon's length in minutes and $x = (x(t))_{0 \le t \le T}$ denote the building's internal temperature in $\,^{\circ}\rm{C}$ during this time.

\subsection{Constraints for the internal temperature}
We suppose that the internal temperature is kept between lower and upper limits $X_{min}$ and $X_{max}$ overnight where $X_{min} < X_{max}$,
\begin{equation}\label{eq:Bulding-Temperature-SLA}
X_{min} \le x(t) \le X_{max},\;\; t \in [0,T].
\end{equation}
A {\it pre-cooling operational strategy} refers to the act of increasing cooling power and using the building's thermal inertia to reduce the need for cooling power at later periods \cite{Reddy1991,Roth2009}. We include pre-cooling in our framework by imposing a constraint on the final temperature value $x(T)$ as follows,
\begin{equation}\label{eq:Bulding-Temperature-Terminal-Constraint}
X_{min} \le x(T) \le \hat{X},
\end{equation}
where $\hat{X}$ in $[X_{min},X_{max}]$ is set by the building manager. Maximum pre-cooling is achieved by setting $\hat{X} = X_{min}$.

\subsection{Calculating the capacity of available decremental reserve}
In this section we outline our methodology for calculating an optimal reference power profile, $C_{ref}$, and an optimal alternative power profile $C_{alt}$. These profiles are used to calculate the optimal instantaneous level of reserve capacity $C_{cap}$ by,
\begin{equation}\label{eq:Reserve-Capacity-Definition}
C_{cap}(t) = C_{alt}(t)-C_{ref}(t),\;\; t \in [0,T].
\end{equation}
If the right-hand side of \eqref{eq:Reserve-Capacity-Definition} is negative then the building is unable to deliver {\it instantaneous} decremental reserve at that time. It is important to note, however, that the building can still deliver {\it total} decremental reserve over the control horizon, provided
\[
\int_{0}^{T}C_{cap}(t){d}t \ge 0.
\]
Negative instantaneous reserve exemplifies a possible consequence of demand response known as the {\em payback effect} \cite{Tindemans2015}, which occurs when the building's cooling equipment recovers after deviating from its normal operation in order to satisfy operational constraints.
\paragraph*{Optimizing the reference profile.}
Suppose there is no request for decremental reserve during $[0,T]$ and the building manager implements the reference profile $C_{ref}$. Letting $P$ ($p / kWh$), where ``p'' stands for pence, denote the positive and constant night-time price of electricity, the total cost to the building manager is,
\begin{equation}\label{eq:Total-Usage-Cost}
\frac{1}{60}\int_{0}^{T}PC_{ref}(t){d}t,
\end{equation}
where we divide by 60 since $T$ is given in minutes. It is reasonable to choose $C_{ref}$ so that it minimizes \eqref{eq:Total-Usage-Cost}, and we formulate an optimal control problem \eqref{Problem:Optimal-Control-Reference-Problem} below to accomplish this.

\paragraph*{Optimizing the alternative profile.}
Let $R$ ($p / kWh$) denote the positive utilization payment received as a reward for the electricity consumed in excess of the reference level $C_{ref}$. When decremental reserve is being delivered according to an alternative profile $C_{alt}$, the {\it instantaneous net cost}
is,
\[
PC_{alt}(t) - R\bigl(C_{alt}(t) - C_{ref}(t)\bigr)^{+},
\]
where $y^{+} = \max(y,0)$. The {\it total net cost} is therefore,
\begin{equation}\label{eq:Total-Net-Cost}
\frac{1}{60}\int_{0}^{T}\bigl[PC_{alt}(t) - R\bigl(C_{alt}(t) - C_{ref}(t)\bigr)^{+}\bigr]{d}t.
\end{equation}
Given the reference profile $C_{ref}$, it is reasonable to choose $C_{alt}$ so that it minimizes the total net cost \eqref{eq:Total-Net-Cost}, and we formulate the second optimal control problem \eqref{Problem:Optimal-Control-Availability-Problem} to achieve this.

\subsection{Delivering the decremental reserve service}
{\it Reserve service instructions} are described by a power profile $C_{ask}$ that indicates how much additional power the building should consume relative to the reference level $C_{ref}(t)$. We suppose $C_{ask}$ is of the form,
\begin{equation}\label{eq:Reserve-Service-Instruction-Profile}
C_{ask}(t) = \sum_{i=1}^{N}C_{i,ask}\mathbf{1}_{[s_{i},e_{i})}(t),
\end{equation}
where $\{(C_{i,ask},s_{i},e_{i})\}_{i = 1,\ldots,N}$ is a sequence of $N \ge 1$ instructions, each consisting of a {\it delivery amount} $C_{i,ask}$ (kW), {\it delivery start time} $s_{i}$ (min), and {\it delivery end time} $e_{i}$ (min), satisfying,
\[
\begin{split}
(i) & \quad 0 \le s_{i} < e_{i} \le T \text{ for } i \in \{1,\ldots,N\},\\
(ii) & \quad e_{i} = s_{i+1} \text{ for } N \ge 2 \text{ and } i \in \{1,\ldots,N-1\},
\end{split}
\]
and $\mathbf{1}_{A}$ denotes the indicator function of a set $A$,
\[
\mathbf{1}_{A}(y) = \begin{cases}
1,& y \in A,\\
0,& y \notin A.
\end{cases}
\]
Therefore, at any time $t \in [s_{i},e_{i})$ the building's power usage must be at least $C_{i,ask}$ (kW) more than the reference level $C_{ref}(t)$. This leads to the following constraint on any delivery profile $C_{del}$ that satisfies the reserve service instructions,
\begin{equation}\label{eq:Reserve-Service-Delivery-Constraint}
C_{del}(t) \ge C_{ref}(t) + C_{i,ask},\;\; s_{i} \le t < e_{i},\;\; i = 1,\ldots,N.
\end{equation}
We assume that the building's instantaneous consumption can be less than the reference level outside of the reserve service instructions' times, allowing it to recover from providing the service as needed. The optimization problem used to determine an optimal $C_{del}$ is formulated as \eqref{Problem:Optimal-Control-Delivery-Problem} below.

\subsection{Internal temperature modelling}\label{Section:Ensemble-Averaging}
\paragraph*{Linear dynamics.}
We assume that $x$ evolves according to the following linear dynamics \cite{Ihara1981,Tindemans2015}:
\begin{equation}
\begin{split}
\dot{x}(t) = -\frac{1}{\tau}\left[x(t) - X_{off} + (X_{off}-X_{on})u(t)\right],\label{eq:Linear-Thermal-Model-Appliance}\\
x(0) \in [X_{on},X_{off}],
\end{split}
\end{equation}
where $\dot{x}$ is the time derivative of $x$ and,
\begin{itemize}
	\item $X_{on}$ and $X_{off}$ are the asymptotic temperatures reached when the cooling equipment operates in the ``on'' and ``off'' states respectively, with $X_{on} < X_{min}$ and $X_{max} < X_{off}$;
	\item $\tau > 0$ is the thermal time constant;
	\item $u(t) \in [0,1]$, the normalized cooling power profile, is the fraction of actual power $C$ used at time $t$,
	\begin{equation}\label{eq:Normalization}
	u(t) = \frac{C(t)}{C_{max}},
	\end{equation}
	where $C_{max}$ (W) is the maximum power usage of the building.
\end{itemize}
Normalized power profiles $u_{ref}$, $u_{alt}$, $u_{del}$, $u_{cap}$, and $u_{ask}$ are defined analogously for the respective reference, alternative, delivery, capacity and reserve instruction profiles.

\begin{definition}
	Let $\mathcal{U}$ denote the set of normalized power profiles $\bigl(u(t)\bigr)_{0 \le t \le T}$ where $u \colon [0,T] \to [0,1]$.
\end{definition}

\paragraph*{Constant temperature control.}
Suppose the temperature $x$ is constant over an interval $[\bar{t}_{0}, \bar{t}_{1}]$ with $0 \le \bar{t}_{0} < \bar{t}_{1} \le T$. Using \eqref{eq:Linear-Thermal-Model-Appliance} any control $\bar{u} \in \mathcal{U}$ that achieves this steady condition satisfies,
\[
-\frac{1}{\tau}\left[x(\bar{t}_{0}) - X_{off} + (X_{off}-X_{on})\bar{u}(t)\right] = 0,\;\;  t \in [\bar{t}_{0},\bar{t}_{1}],
\]
and therefore,
\begin{equation}\label{eq:Steady-State-Control}
\bar{u}(t) = \bar{u}(\bar{t}_{0}) = \frac{X_{off} - x(\bar{t}_{0})}{X_{off}-X_{on}},\;\;  t \in [\bar{t}_{0},\bar{t}_{1}].
\end{equation}

\paragraph*{Analytic solution for step controls.}
Suppose $u \in \mathcal{U}$ satisfies, \begin{equation}\label{eq:Piecewise-Constant-Control}
u(t) = \sum_{k = 1}^{n_{p}}u_{k}\mathbf{1}_{[t_{k-1},t_{k})}(t),\;\; t \in [0,T],
\end{equation}
where $n_{p} > 1$ is an integer and $\{t_{k}\}_{k = 0}^{n_{p}}$ is a sequence of time points $0 = t_{0} < \ldots < t_{n_{p}} = T$ that partition the control horizon $[0,T]$ into $n_{p} > 0$ contiguous subintervals. The solution $(x(t))_{0 \le t \le T}$ to \eqref{eq:Linear-Thermal-Model-Appliance} corresponding to the control \eqref{eq:Piecewise-Constant-Control} is continuous and satisfies,
\begin{equation}\label{eq:Analytic-Solution-Piecewise-Constant-Control}
\begin{split}
x(t) = e^{-\frac{t - t_{k-1}}{\tau}}x(t_{k-1}) + \Bigl(1 - e^{-\frac{t-t_{k-1}}{\tau}}\Bigr)\bigl(X_{off} + (X_{on} - X_{off})u_{k}\bigr),\\
t_{k-1} \le t < t_{k},\;\; 1 \le k \le n_{p}.
\end{split}
\end{equation}

\subsection{Three optimal control problems for off-peak decremental replacement reserve provision}

\paragraph*{Problem 1. Optimal reference power usage.}
The reference power profile $C_{ref}$ should conduce minimum expenditure if decremental reserve is not requested over the control horizon. The following optimal control problem, formulated using the normalized reference profile, is suitable for our application.
\begin{gather}
\text{minimize}\; \int_{0}^{T}\bigl[u_{ref}(t) + \alpha_{ref}\bigl(u_{ref}(t)\bigr)^{2}\bigr]{d}t\;\; \text{over}\;\; u_{ref} \in \mathcal{U}\;\; \text{subject to:} \nonumber \\
\begin{split}
(i) & \quad \dot{x}(t) = f(t,x(t),u_{ref}(t)) \text{ given by } \eqref{eq:Linear-Thermal-Model-Appliance}, \\
(ii) & \quad X_{min} \le x(t) \le X_{max},\;\; t \in [0,T],\\
(iii) & \quad x(0) \in [X_{min},X_{max}] \text{ and } x(T) \in [X_{min},\hat{X}], \\
\end{split}\label{Problem:Optimal-Control-Reference-Problem}
\end{gather}
where $\alpha_{ref} > 0$ is a constant that weighs the importance of the {\it regularization term} $\bigl(u_{ref}(t)\bigr)^{2}$. The regularizer is used in \eqref{Problem:Optimal-Control-Reference-Problem} to disfavour solutions where $u_{ref}$ alternates rapidly between its minimum and maximum possible values. In the absence of this regularizer, theory states that this unwanted behaviour can be optimal in \eqref{Problem:Optimal-Control-Reference-Problem} since the control variable $u_{ref}$ then appears linearly in both the objective function and state dynamics \cite{Maurer1977}.

\paragraph*{Problem 2. Optimal alternative power usage.}
Given the normalized reference power profile $u_{ref}$, night-time electricity price $P$, and utilization payment $R$,
\begin{gather}
\text{minimize}\; J(u_{alt}\,;\,u_{ref},P,R)\;\; \text{over}\;\; u_{alt} \in \mathcal{U}\;\; \text{subject to:} \nonumber \\
\begin{split}
(i) & \quad \dot{x}(t) = f(t,x(t),u_{alt}(t)) \text{ given by } \eqref{eq:Linear-Thermal-Model-Appliance}, \\
(ii) & \quad X_{min} \le x(t) \le X_{max},\;\; t \in [0,T],\\
(iii) & \quad x(0) \in [X_{min},X_{max}] \text{ and } x(T) \in [X_{min},\hat{X}], \\
(iv) & \quad \int_{0}^{T}\bigl[u_{alt}(t) -  u_{ref}(t) - \tfrac{R}{P}\bigl(u_{alt}(t) - u_{ref}(t)\bigr)^{+}\bigr]{d}t \le 0, \\
\end{split}\label{Problem:Optimal-Control-Availability-Problem}
\end{gather}
where $J(u_{alt}\,;\,u_{ref},P,R)$ is given by,
\begin{equation}\label{eq:Optimal-Control-Availability-Problem-Objective}
\begin{split}
J(u_{alt}\,;\,u_{ref},P,R) = {} & \int_{0}^{T}\bigl[u_{alt}(t) - \tfrac{R}{P}\bigl(u_{alt}(t) - u_{ref}(t)\bigr)^{+}\bigr]{d}t \\
& + \alpha_{alt}\int_{0}^{T}\bigl(u_{alt}(t)\bigr)^{2}{d}t,
\end{split}
\end{equation}
and, similar to \eqref{Problem:Optimal-Control-Reference-Problem} above, $\alpha_{alt} > 0$ is a constant used to weigh the importance of a quadratic regularizer. The effect this parameter can have on the results is illustrated in the Appendix. Constraint \eqref{Problem:Optimal-Control-Availability-Problem}-(iv) ensures that the building manager is not worse off financially by following $u_{alt}$ instead of $u_{ref}$, and is always satisfied when $R \ge P$. If ``$0$'' on the right-hand side of \eqref{Problem:Optimal-Control-Availability-Problem}--(iv) is replaced by ``$-\gamma$'' where $\gamma \ge 0$ then a solution to \eqref{Problem:Optimal-Control-Availability-Problem} guarantees the building manager a minimum total net profit of $\gamma\bigl(\frac{C_{max}P}{60}\bigr)$ pence relative to the reference profile's cost.

\paragraph*{Problem 3. Optimal reserve service delivery.}
Since the building manager is only compensated for the additional demand as instructed, it is reasonable to require that the building uses no more power than that needed to satisfy the reserve instructions and internal temperature limits. Using the normalized profiles $u_{ref}$ and $u_{ask}$, we therefore formulate the reserve service delivery problem as follows.
\begin{gather}
\text{minimize}\; \int_{0}^{T}\bigl[u_{del}(t) + \alpha_{del}\bigl(u_{del}(t)\bigr)^{2}\bigr]{d}t\;\; \text{over}\;\; u_{del} \in \mathcal{U}\;\; \text{subject to:} \nonumber \\
\begin{split}
(i) & \quad \dot{x}(t) = f(t,x(t),u_{del}(t)) \text{ given by } \eqref{eq:Linear-Thermal-Model-Appliance},\\
(ii) & \quad X_{min} \le x(t) \le X_{max},\;\; t \in [0,T],\\
(iii) & \quad x(0) \in [X_{min},X_{max}] \text{ and } x(T) \in [X_{min},\hat{X}],\\
(iv) & \quad u_{del}(t) \ge u_{ref}(t) + u_{i,ask},\;\; s_{i} \le t < e_{i},\;\; i = 1,\ldots,N,\\
\end{split}\label{Problem:Optimal-Control-Delivery-Problem}
\end{gather}
where, similar to \eqref{Problem:Optimal-Control-Reference-Problem} above, $\alpha_{del} > 0$ is a constant used to weigh the importance of the quadratic regularizer. For computations and illustrations, it is convenient to represent the constraint \eqref{Problem:Optimal-Control-Delivery-Problem}--(iv) by its equivalent form,
\[
u_{del}(t) \ge u_{ins}(t),\;\; t \in [0,T],
\]
where $u_{ins}$ is the profile of minimal instructed power usage,
\[
u_{ins}(t) = \sum_{i=1}^{N}\bigl[u_{ref}(t) + u_{i,ask}\bigr]\mathbf{1}_{[s_{i},e_{i})}(t),
\]
and we used the property $u_{del}(t) \ge 0$.

\section{Numerical simulations}\label{Section:Numerical-Results}
Theory, for example \cite{Cesari1983,Clarke2013}, guarantees the existence of a solution to each of the problems \eqref{Problem:Optimal-Control-Reference-Problem}, \eqref{Problem:Optimal-Control-Availability-Problem}, and \eqref{Problem:Optimal-Control-Delivery-Problem}. In this section we present results of the associated numerical solutions, which we obtained using the control parametrization method \cite{Teo1991,Rehbockt1999} outlined in Appendix~\ref{Section:Control-Parametrization}. The hypothetical building for the experiments has computing equipment that generates a significant amount of electricity demand for temperature cooling. The control horizon is $T = 360$ minutes long, lasting from 00:00 to 06:00. We assume an internal temperature range of $X_{min} = 18$ ($^{\circ}\rm{C}$) to $X_{max} = 27$ ($^{\circ}\rm{C}$) must be maintained during this time, which is the recommended range in the ASHRAE Standard 90.4-2016 for data centres \cite{ASHRAE2016}. For the simulations we set the utilization payment at either 75\%, 100\% or 125\% of the electricity price. According to the 2017 DTU market information \cite{DemandTurnUp}, the utilization payment was typically in the range $R \in [6, 10]$ (p / kWh). Table \ref{Table:Params} lists values for the parameters used to generate the numerical results.

\begin{table}[!ht]
	\renewcommand{\arraystretch}{1.3}
	\caption{Parameters for numerical simulations.}
	\label{Table:Params}
	\centering
	\begin{tabular}{ccc}
		\bf Parameter & \bf Value & \bf Units \\
		\hline 
		$T$ & $360$ & minutes \\
		$\tau$ & $120$ & minutes \\
		$X_{off}$ & $35$ & $\,^{\circ}\rm{C}$ \\
		$X_{on}$ & $10$ & $\,^{\circ}\rm{C}$ \\
		$x(0)$ & $27$ & $\,^{\circ}\rm{C}$ \\
		$X_{min}$ & $18$ & $\,^{\circ}\rm{C}$ \\
		$X_{max}$ & $27$ & $\,^{\circ}\rm{C}$ \\
		$\frac{R}{P}$ & $\frac{R}{P} \in \{\frac{3}{4},1,\frac{5}{4}\}$ & $1$ \\
		\hline
	\end{tabular}
\end{table}

\subsection{Optimal reference power usage}
In this section we highlight important characteristics of the control and temperature trajectories corresponding to an optimal reference power profile $u_{ref}$. These trajectories are intuitive and show that the optimal control uses minimal effort to keep the temperature within its constraints. In particular,
\begin{enumerate}
	\item Starting at the maximum feasible level $X_{max}$, the temperature is kept constant at this level until a time $\bar{t}_{1} \in (0,T)$;
	\item From time $\bar{t}_{1}$ onwards, power consumption is increased steadily until it reaches the maximum level at a time $\bar{t}_{2} \in (\bar{t}_{1},T)$;
	\item From time $\bar{t}_{2}$ to $T$, maximum power is used to steer the temperature to $\hat{X}$, the upper limit of feasible values at time $T$.
\end{enumerate}
Consequently, when following the optimal reference profile $u_{ref}$, the building has capacity to provide decremental reserve on $[0,\bar{t}_{1}]$, but is unavailable to provide this service from time $\bar{t}_{2}$ onwards when reference power usage is at its highest. Note that these optimal characteristics are sensible for temperature dynamics other than the linear one \eqref{eq:Linear-Thermal-Model-Appliance}.

The linear model allows for easy approximation of the terms described above. For example, using \eqref{eq:Steady-State-Control} the control $u_{ref}$ satisfies,
\begin{equation}\label{eq:Constant-Optimal-Control}
u_{ref}(t) \approx \frac{X_{off} - X_{max}}{X_{off}-X_{on}},\;\;  t \in [0,\bar{t}_{1}].
\end{equation}
Using \eqref{eq:Analytic-Solution-Piecewise-Constant-Control}, the time $\bar{t}_{2}$ in the description approximately satisfies,
\[
\hat{X} = \bigl(e^{-\frac{T - \bar{t}_{2}}{\tau}}\bigr)X_{max} + \bigl(1 - e^{-\frac{T-\bar{t}_{2}}{\tau}}\bigr)X_{on},
\]
which we solve to get,
\begin{equation}\label{eq:Approximate-Full-Power-Time}
\bar{t}_{2} \approx T - \tau\log\left(\frac{X_{max} - X_{on}}{\hat{X} - X_{on}}\right).
\end{equation}
The time $\bar{t}_{1}$ can be determined by calculating how long it takes to go from the constant power level at $\bar{t}_{1}$ to full power at $\bar{t}_{2}$ using \eqref{eq:Constant-Optimal-Control} and \eqref{eq:Approximate-Full-Power-Time}. Using \eqref{eq:Approximate-Full-Power-Time} we see that the duration of unavailability, $T - \bar{t}_{2}$, is proportional to the building's thermal time constant $\tau$, and also increases with the level of pre-cooling at time $T$, which is controlled by the parameter $\hat{X}$.

\subsection{Optimal capacity of available reserve}
\begin{figure*}[!ht]
	\centering
	\subfigure[$\frac{R}{P} = \frac{3}{4}$]{%
		\resizebox*{6.15cm}{!}{\includegraphics{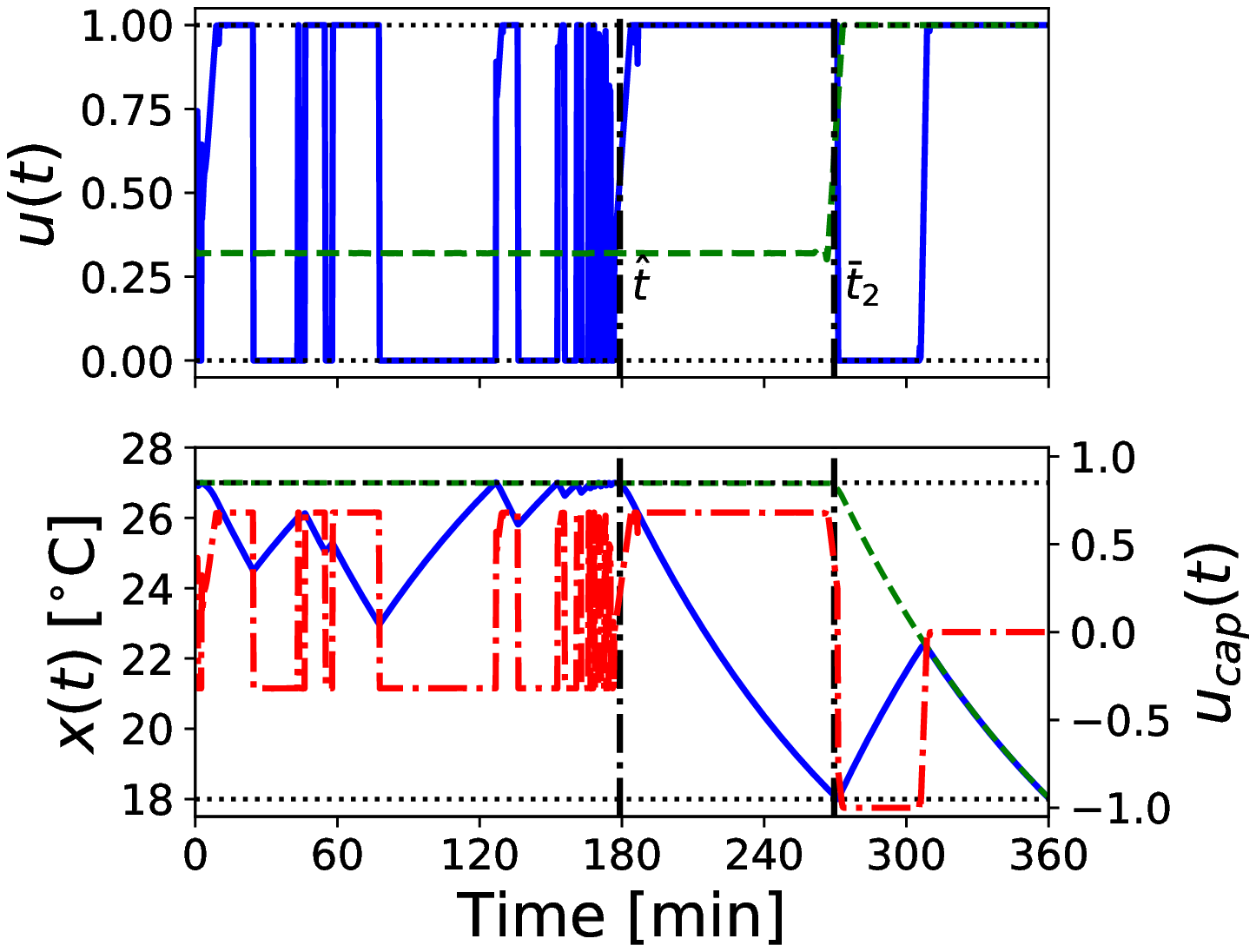}}
		\label{fig:Optimized-Alternative-Power-1}
	}
	\hspace{5pt}
	\subfigure[$\frac{R}{P} = 1$]{%
		\resizebox*{6cm}{!}{\includegraphics{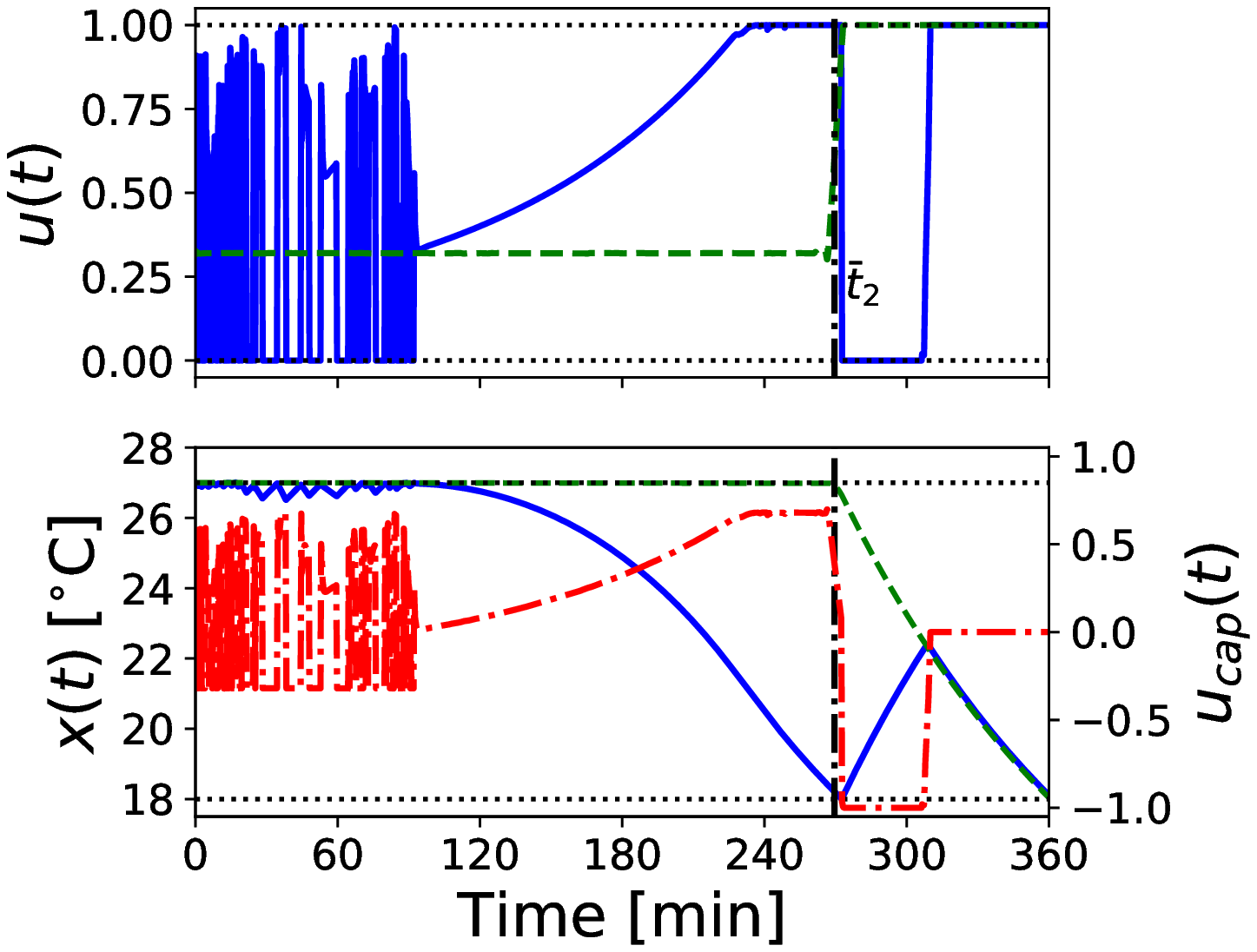}}
		\label{fig:Optimized-Alternative-Power-2}
	}\hspace{5pt}
	\subfigure[$\frac{R}{P} = \frac{5}{4}$]{%
		\resizebox*{6.5cm}{!}{\includegraphics{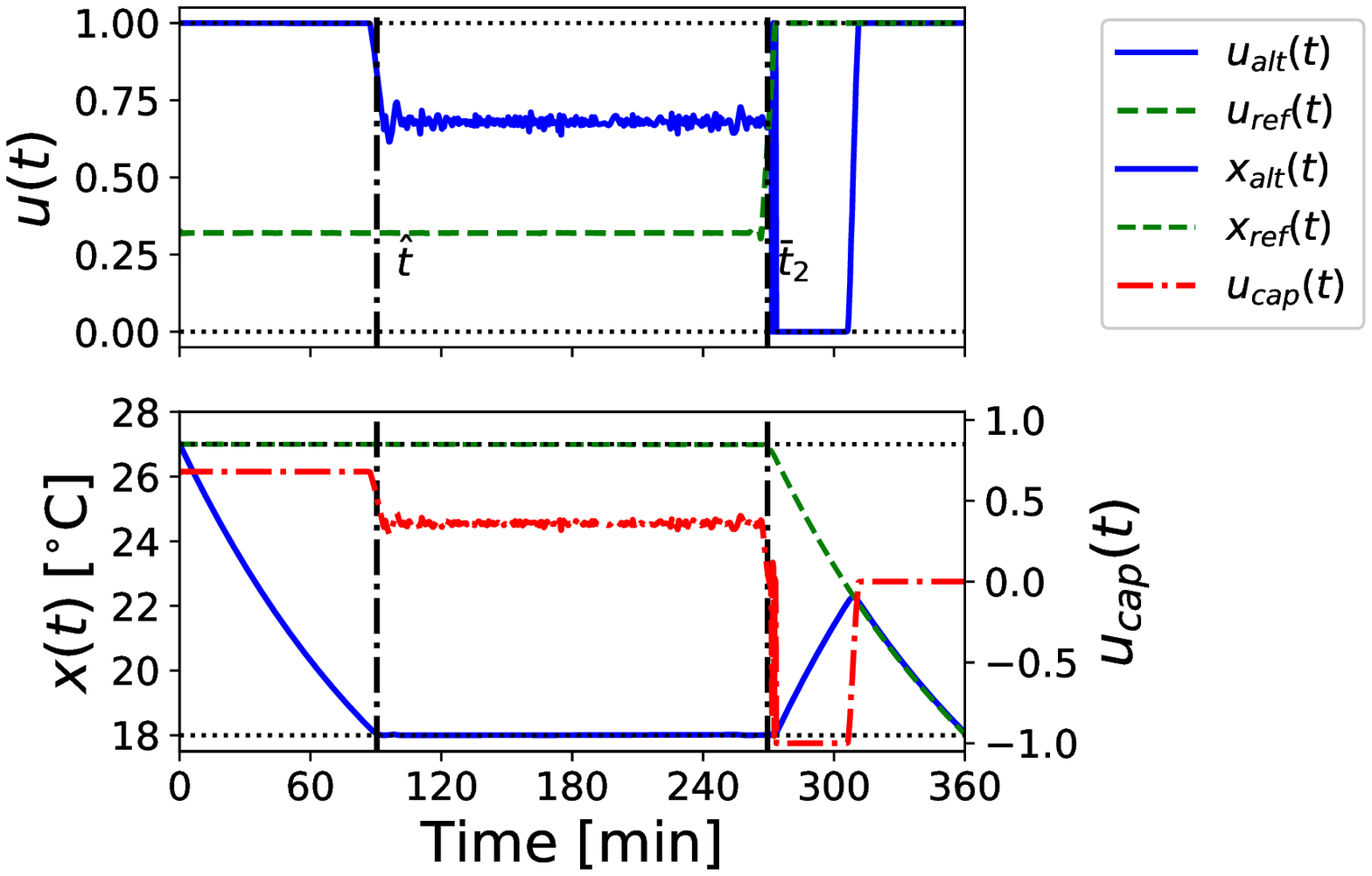}}
		\label{fig:Optimized-Alternative-Power-3}
	}\hspace{5pt}
	\subfigure[Normalized net profit]{%
		\resizebox*{5.5cm}{!}{\includegraphics[height=0.225\textheight, width=0.475\textwidth]{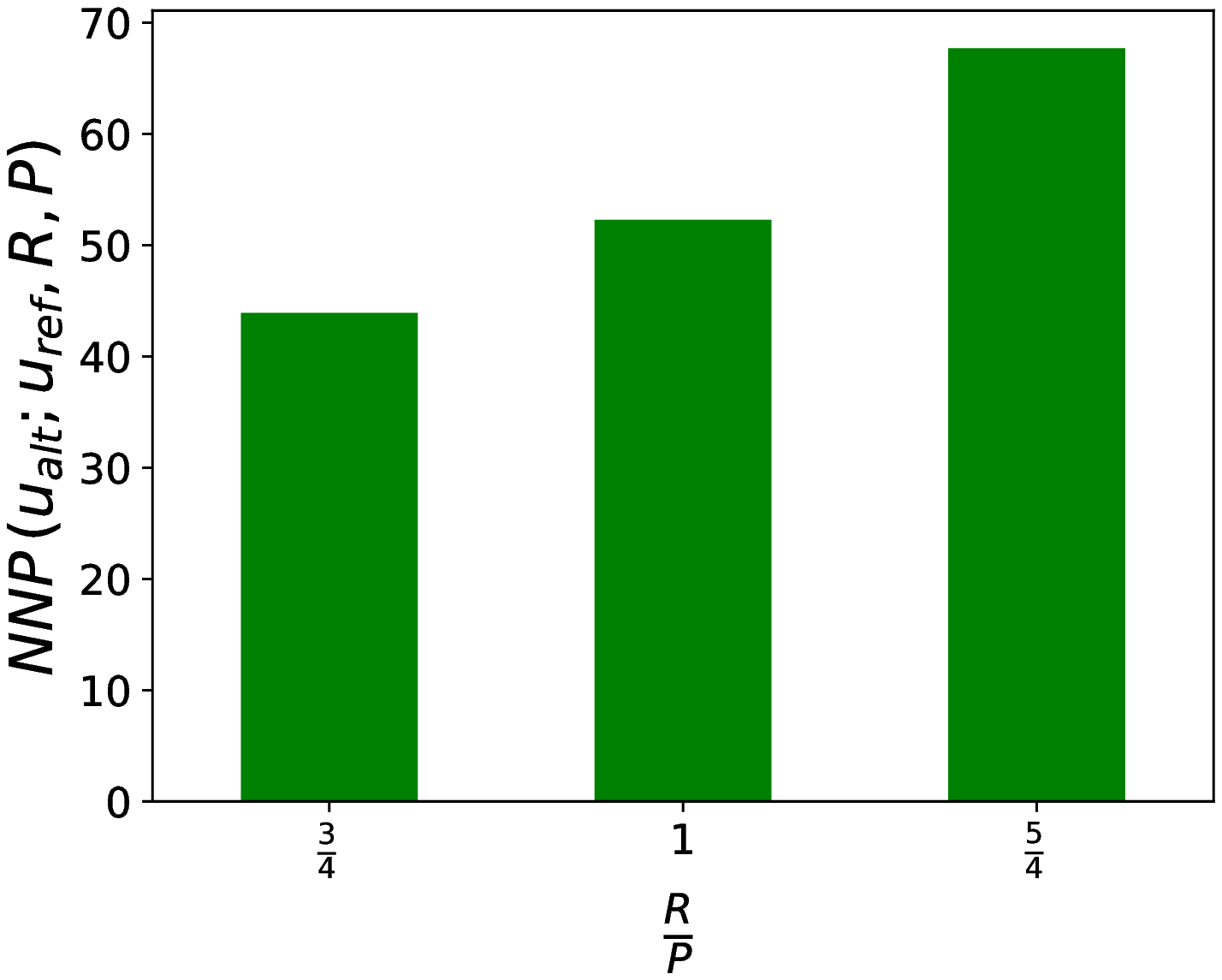}}
		\label{fig:Optimized-Alternative-Power-Net-Profits}
	}
	\caption{Optimized internal temperature $x$, alternative power usage $u_{alt}$ and reserve capacity $u_{cap}$ corresponding to three different values for the benefit-cost ratio $\frac{R}{P}$. The time $\bar{t}_{2}$ (cf. \eqref{eq:Approximate-Full-Power-Time}) at which the reference profile starts to apply maximum power near the control horizon's end is shown in each case. In (a) $\frac{R}{P} = \frac{3}{4}$ and the building tends to be available for decremental reserve at maximum capacity on short intervals. The longest duration of reserve occurs on an interval $[\hat{t},\bar{t}_{2}]$ where the temperature is cooled from the maximum allowed value to the minimum one. In (b) $\frac{R}{P} = 1$ and the building's capacity for decremental reserve is gently increased up to the maximum value before $\bar{t}_{2}$. In (c) $\frac{R}{P} = \frac{5}{4}$ and the building is mostly available for decremental reserve. On average, less than maximum power is sustained on an interval $[\check{t},\bar{t}_{2}]$ while the temperature is at its minimum allowed value. Figure (d) shows that, as expected, the normalized net profit increases with $\frac{R}{P}$. Parameter values are $x(0) = 28$ for the initial temperature, $\hat{X} = 18$ for the pre-cooling value, and $\alpha_{ref} = \alpha_{alt} = 0.01$ for the control regularizers. Temperature and control constraints are shown using dotted horizontal lines.}
	\label{fig:Optimized-Alternative-Power}
\end{figure*}

This section illustrates the numerical results for the reserve capacity problem \eqref{Problem:Optimal-Control-Availability-Problem}. Figure~\ref{fig:Optimized-Alternative-Power} shows the alternative power usage and capacity profiles corresponding to three different values of $\frac{R}{P}$. It also shows the normalized net profit relative to the reference profile for the three cases, which is defined as the negative of constraint \eqref{Problem:Optimal-Control-Availability-Problem}-(iv),
\[
NNP(u_{alt};u_{ref},R,P) = \int_{0}^{T}\bigl[u_{ref}(t) - u_{alt}(t) + \tfrac{R}{P}\bigl(u_{alt}(t) - u_{ref}(t)\bigr)^{+}\bigr]{d}t.
\]
Multiplying this value by $\bigl(\frac{C_{max}P}{60}\bigr)$ gives the total net profit relative to the reference profile in pence.

\subsubsection*{Case 1: $\frac{R}{P} = \frac{3}{4}$}

In Figure~\ref{fig:Optimized-Alternative-Power-1} the benefit-cost ratio satisfies $\frac{R}{P} < 1$, and the optimal control tends to use short bursts of pre-cooling to minimize the cost of providing reserve capacity. Consequently, there are many intervals during which the building is either unavailable for decremental reserve, or available for short periods at maximum capacity. The longest period of sustained maximal capacity occurs just before time $\bar{t}_{2}$ when the reference profile power usage is at its highest. The normalized level of sustained capacity during this time is approximately,
\begin{equation}\label{eq:Case-1-Normalized-Level-Sustained-Capacity}
u_{cap}(t) \approx \frac{X_{max} - X_{on}}{X_{off}-X_{on}},\;\;  t \in [\hat{t},\bar{t}_{2}],
\end{equation}
where, analogous to \eqref{eq:Approximate-Full-Power-Time} above for the reference profile, $\hat{t}$ is given by,
\begin{equation}\label{eq:Case-1-Approximate-Full-Power-Time}
\hat{t} \approx \bar{t}_{2} - \tau\log\left(\frac{X_{max} - X_{on}}{X_{min} - X_{on}}\right).
\end{equation}

Using the expression for $\bar{t}_{2}$ in \eqref{eq:Approximate-Full-Power-Time} and $\hat{X} = X_{min}$ shows,
\[
\bar{t}_{2} - \hat{t} \approx  T - \bar{t}_{2},
\]
and $u_{alt}$ essentially time-shifts the final period of maximum consumption that occurred under the reference profile $u_{ref}$. It is important to note that the net cost of providing maximum reserve capacity increases as $R$ decreases. In the presence of thermal losses, this means there may no longer be any profitable solutions to \eqref{Problem:Optimal-Control-Availability-Problem} of this form if $R$ is too low.

\subsubsection*{Case 2: $R = P$}
In Figure \ref{fig:Optimized-Alternative-Power-2} the benefit-cost ratio satisfies $\frac{R}{P} = 1$ and the alternative profile $u_{alt}$ exhibits complex bang-bang behaviour early within the control horizon. This complex behaviour can be explained by noticing that when $\frac{R}{P} = 1$ the instantaneous cost in \eqref{eq:Optimal-Control-Availability-Problem-Objective} satisfies,
\begin{equation}\label{eq:When-Reward-Equals-Cost}
\begin{split}
u_{alt}(t) - \tfrac{R}{P}\bigl(u_{alt}(t) - u_{ref}(t)\bigr)^{+} + \alpha_{alt}\bigl(u_{alt}(t)\bigr)^{2} = {} & \min\bigl(u_{alt}(t),u_{ref}(t)\bigr) \\
& + \alpha_{alt} \bigl(u_{alt}(t)\bigr)^{2}.
\end{split}
\end{equation}
When $\alpha_{alt}$ is low, the term $\min\bigl(u_{alt}(t),u_{ref}(t)\bigr)$ can dominate \eqref{eq:When-Reward-Equals-Cost} and incentivize $u_{alt}$ to use bang-bang control behaviour in order to be more cost effective than the reference profile $u_{ref}$, which already uses minimal effort. This period of complex control behaviour is clearly unsuitable for providing decremental reserve. Moreover, even after this period the profile $u_{alt}$ may not be practical since, unlike in Figure~\ref{fig:Optimized-Alternative-Power-1}, it {\em slowly} ramps up to maximum consumption and provides only a short period of constant reserve before $\bar{t}_{2}$.

\subsubsection*{Case 3: $\frac{R}{P} = \frac{5}{4}$}

In Figure~\ref{fig:Optimized-Alternative-Power-3} the benefit-cost ratio satisfies $\frac{R}{P} > 1$. The control $u_{alt}$ initially applies maximum power to steer the internal temperature to its minimum allowable value $X_{min}$. Let $\check{t}$ denote the first time that the temperature hits $X_{min}$. Analogous to \eqref{eq:Case-1-Approximate-Full-Power-Time}, $\check{t}$ is given approximately by,
\begin{equation}\label{eq:Case-2-Approximate-Full-Power-Time}
\check{t} \approx \tau\log\left(\frac{X_{max} - X_{on}}{X_{min} - X_{on}}\right),
\end{equation}
which we recall is equal to $T - \bar{t}_{2}$, the length of time that the reference profile $u_{ref}$ applies maximum power. The control $u_{alt}$ keeps the temperature at $X_{min}$ from $\check{t}$ until the time $\bar{t}_{2}$ by applying power that is on average equal to (cf. \eqref{eq:Steady-State-Control}),
\[
u_{alt}(t) \approx \frac{X_{off} - X_{min}}{X_{off}-X_{on}},\;\;  t \in [\check{t},\bar{t}_{2}],
\]
and this control provides a corresponding sustained level of capacity that is approximately,
\begin{equation}\label{eq:Case-2-Normalized-Level-Sustained-Capacity}
u_{cap}(t) \approx \frac{X_{max} - X_{min}}{X_{off}-X_{on}},\;\;  t \in [\check{t},\bar{t}_{2}].
\end{equation}

\subsection{Optimal reserve service delivery}

In this section we summarize the numerical results for the optimal reserve service delivery problem \eqref{Problem:Optimal-Control-Delivery-Problem}. For further details see Appendix \ref{Appendix:Trajectories}. The optimized delivery profile uses minimal power to satisfy the temperature constraints and reserve service instructions. After delivering the service, the internal temperature is lower than the reference level, and the cooling equipment can be turned off temporarily while the temperature rises within its permitted range. The building manager therefore benefits doubly, by receiving the utilization payment for delivering the reserve service, and by reducing the electricity usage consequent to pre-cooling. The cooling equipment remains off for a longer duration if less pre-cooling is required at the control horizon's end.

\section{Summary and recommendations}\label{Section:Conclusion}

As the energy transition transforms power grids across the globe, high levels of intermittent renewable generation complicate the job of continuously balancing power supply and demand, which is necessary for the grid's stability. New ancillary services have emerged in this regard, such as National Grid UK's Demand Turn Up (DTU) \cite{DemandTurnUp}, which is a reserve service that incentivizes large energy consumers to increase their electricity demand, for example during overnight periods of high output from renewable generation and low overall demand. In this paper we explore the optimal participation of a commercial building, through the control of its temperature cooling equipment, in such an ancillary service initiative. We provide a computational framework for solving this problem that takes into account the economic incentives given. The framework has three main outputs:
\begin{enumerate}
	\item an optimal reference night-time control schedule for the cooling system when it does not provide DTU.
	\item an optimal schedule of DTU capacity relative to the reference for a given remuneration.
	\item an optimal night-time control schedule to fulfil DTU instructions.
\end{enumerate}
The framework also takes into account the building's relaxation dynamics, so that DTU requests are used as an opportunity to optimally pre-cool the building. In addition to the DTU payment, this pre-cooling reduces energy usage during the subsequent morning peak period, which is a financial benefit to the customer and also reduces stress on the grid.

The optimal control schedule used as reference or to fulfil the DTU instructions is intuitive and uses minimal effort to satisfy the temperature and power constraints. Consistent with studies such as \cite{DeConinck2016, Blum2017}, we find that the level of participation in the ancillary service is affected by the dynamics and constraints for the internal temperature, and how well the building manager is remunerated. Figure \ref{fig:Optimized-Alternative-Power} shows that the building's capacity for DTU depends crucially on the {\it benefit-cost ratio} $\frac{R}{P}$, where $R$ is the utilization payment for DTU and $P$ is the night-time price of electricity.
\begin{itemize}
	\item When $\frac{R}{P} < 1$ the optimal alternative control schedule has a complicated structure with time-shifts in consumption causing intermittent periods of constant reserve throughout the control horizon. Nevertheless, this schedule can be partially implemented as there is a sufficiently long period of constant reserve before the control horizon's end.
	\item When $\frac{R}{P} = 1$ the optimal alternative control schedule also has a complicated structure but is not implementable since it lacks sufficiently long periods of constant reserve.
	\item When $\frac{R}{P} > 1$ the optimal alternative control schedule has the least complicated structure and sustains constant reserve for long periods throughout the control horizon. It can consist entirely of two contiguous intervals of constant reserve, which is consistent with what is allowed in practice \cite{DemandTurnUp}.
\end{itemize}
Even in Case 1 above, where the utilization payment is lower than the night-time electricity price ($\frac{R}{P} = \frac{3}{4}$), participation in DTU was found to be profitable. This is because simply shifting pre-planned HVAC operation to a time earlier in the night leads to increased demand at the earlier time, attracting compensation under DTU. However when $\frac{R}{P} < 1$ the optimal control strategy becomes more complex as $R$ decreases, fluctuating more frequently between minimum and maximum power as it ``hunts'' for profit. Demand reductions due to payback effects are undesirable as they undermine the purpose of DTU, and frequent rapid power fluctuations may be problematic for system stability. Therefore, in order to economically incentivize a building manager to provide DTU in a practical way, the level of remuneration must be sufficiently high.

Possible future extensions of our work would include controlling an ensemble of possibly heterogeneous thermostatic loads \cite{Tindemans2015} and considering measured temperature relaxation dynamics for each member of the ensemble. The heterogeneity can come from different rooms in a single building, or from an aggregation of multiple buildings. Our model may also be extended to consider the uncertainty in parameters affecting the internal temperature, such as the external weather conditions, or the uncertainty in being called to provide DTU.

\begin{thebibliography}{10}
	\providecommand{\url}[1]{#1}
	\csname url@samestyle\endcsname
	\providecommand{\newblock}{\relax}
	\providecommand{\bibinfo}[2]{#2}
	\providecommand{\BIBentrySTDinterwordspacing}{\spaceskip=0pt\relax}
	\providecommand{\BIBentryALTinterwordstretchfactor}{4}
	\providecommand{\BIBentryALTinterwordspacing}{\spaceskip=\fontdimen2\font plus
		\BIBentryALTinterwordstretchfactor\fontdimen3\font minus
		\fontdimen4\font\relax}
	\providecommand{\BIBforeignlanguage}[2]{{%
			\expandafter\ifx\csname l@#1\endcsname\relax
			\typeout{** WARNING: IEEEtran.bst: No hyphenation pattern has been}%
			\typeout{** loaded for the language `#1'. Using the pattern for}%
			\typeout{** the default language instead.}%
			\else
			\language=\csname l@#1\endcsname
			\fi
			#2}}
	\providecommand{\BIBdecl}{\relax}
	\BIBdecl
	
	\bibitem{Kirschen2004}
	D.~S. Kirschen and G.~Strbac, \emph{{Fundamentals of Power System
			Economics}}.\hskip 1em plus 0.5em minus 0.4em\relax Chicester, England: John
	Wiley {\&} Sons, 2004.
	
	\bibitem{Szabo2017}
	D.~Z. Szab{\'{o}} and R.~Martyr, ``{Real option valuation of a decremental
		regulation service provided by electricity storage},'' \emph{Philosophical
		Transactions of the Royal Society A: Mathematical, Physical and Engineering
		Sciences}, vol. 375, no. 2100, aug 2017.
	
	\bibitem{Rothleder2014}
	M.~Rothleder and C.~Loutan, ``{Case Study - Renewable Integration: Flexibility
		Requirement, Potential Overgeneration, and Frequency Response Challenges},''
	in \emph{Renewable Energy Integration - Practical Management of Variability,
		Uncertainty, and Flexibility in Power Grids}, L.~E. Jones, Ed.\hskip 1em plus
	0.5em minus 0.4em\relax New York: Elsevier Inc., 2014, ch.~6, pp. 67--79.
	
	\bibitem{Hirst1998}
	E.~Hirst and B.~Kirby, ``{Unbundling Generation and Transmission Services for
		Competitive Electricity Markets: Examining Ancillary Services},'' The
	National Regulatory Research Institute, The Ohio State University, Columbus,
	OH, USA, Tech. Rep., jan 1998.
	
	\bibitem{Vardakas2015}
	J.~S. Vardakas, N.~Zorba, and C.~V. Verikoukis, ``{A Survey on Demand Response
		Programs in Smart Grids: Pricing Methods and Optimization Algorithms},''
	\emph{IEEE Communications Surveys {\&} Tutorials}, vol.~17, no.~1, pp.
	152--178, 2015.
	
	\bibitem{Paterakis2017}
	N.~G. Paterakis, O.~Erdin{\c{c}}, and J.~P. Catal{\~{a}}o, ``{An overview of
		Demand Response: Key-elements and international experience},''
	\emph{Renewable and Sustainable Energy Reviews}, vol.~69, pp. 871--891, mar
	2017.
	
	\bibitem{Olivieri2014}
	S.~J. Olivieri, G.~P. Henze, C.~D. Corbin, and M.~J. Brandemuehl, ``{Evaluation
		of commercial building demand response potential using optimal short-term
		curtailment of heating, ventilation, and air-conditioning loads},''
	\emph{Journal of Building Performance Simulation}, vol.~7, no.~2, pp.
	100--118, mar 2014.
	
	\bibitem{Pavlak2014}
	G.~S. Pavlak, G.~P. Henze, and V.~J. Cushing, ``{Optimizing commercial building
		participation in energy and ancillary service markets},'' \emph{Energy and
		Buildings}, vol.~81, pp. 115--126, oct 2014.
	
	\bibitem{Lawrence2016}
	T.~M. Lawrence, M.-C. Boudreau, L.~Helsen, G.~Henze, J.~Mohammadpour,
	D.~Noonan, D.~Patteeuw, S.~Pless, and R.~T. Watson, ``{Ten questions
		concerning integrating smart buildings into the smart grid},'' \emph{Building
		and Environment}, vol. 108, pp. 273--283, nov 2016.
	
	\bibitem{Kim2016}
	Y.-J. Kim, D.~H. Blum, N.~Xu, L.~Su, and L.~K. Norford, ``{Technologies and
		Magnitude of Ancillary Services Provided by Commercial Buildings},''
	\emph{Proceedings of the IEEE}, vol. 104, no.~4, pp. 758--779, apr 2016.
	
	\bibitem{DeConinck2016}
	R.~{De Coninck} and L.~Helsen, ``{Quantification of flexibility in buildings by
		cost curves – Methodology and application},'' \emph{Applied Energy}, vol.
	162, pp. 653--665, jan 2016.
	
	\bibitem{Blum2017}
	D.~H. Blum, T.~Zakula, and L.~K. Norford, ``{Opportunity Cost Quantification
		for Ancillary Services Provided by Heating, Ventilating, and Air-Conditioning
		Systems},'' \emph{IEEE Transactions on Smart Grid}, vol.~8, no.~3, pp.
	1264--1273, may 2017.
	
	\bibitem{Dayarathna2016}
	M.~Dayarathna, Y.~Wen, and R.~Fan, ``{Data Center Energy Consumption Modeling:
		A Survey},'' \emph{IEEE Communications Surveys {\&} Tutorials}, vol.~18,
	no.~1, pp. 732--794, 2016.
	
	\bibitem{DemandTurnUp}
	\BIBentryALTinterwordspacing
	{National Grid UK}, ``{Demand Turn Up},'' 2018. [Online]. Available:
	\url{https://www.nationalgrid.com/uk/electricity/balancing-services/reserve-services/demand-turn}
	\BIBentrySTDinterwordspacing
	
	\bibitem{Clarke2013}
	F.~Clarke, \emph{{Functional Analysis, Calculus of Variations and Optimal
			Control}}, ser. Graduate Texts in Mathematics.\hskip 1em plus 0.5em minus
	0.4em\relax London, England: Springer London, 2013, vol. 264.
	
	\bibitem{Teo1991}
	K.~L. Teo and C.~J. Goh, \emph{{A unified computational approach to optimal
			control problems}}.\hskip 1em plus 0.5em minus 0.4em\relax Essex, England:
	Longman Scientic {\&} Technical, 1991.
	
	\bibitem{Rehbockt1999}
	V.~Rehbockt, K.~L. Teo, L.~S. Jennings, and H.~W.~J. Lee, ``{A Survey of the
		Control Parametrization and Control Parametrization Enhancing Methods for
		Constrained Optimal Control Problems},'' in \emph{Progress in Optimization
		(Applied Optimization Vol 30)}.\hskip 1em plus 0.5em minus 0.4em\relax
	Springer US, 1999, pp. 247--275.
	
	\bibitem{Deng2015}
	R.~Deng, Z.~Yang, M.-Y. Chow, and J.~Chen, ``{A Survey on Demand Response in
		Smart Grids: Mathematical Models and Approaches},'' \emph{IEEE Transactions
		on Industrial Informatics}, vol.~11, no.~3, pp. 570--582, jun 2015.
	
	\bibitem{Wang2008}
	S.~Wang and Z.~Ma, ``{Supervisory and Optimal Control of Building HVAC Systems:
		A Review},'' \emph{HVAC{\&}R Research}, vol.~14, no.~1, pp. 3--32, jan 2008.
	
	\bibitem{Reddy1991}
	T.~Reddy, L.~Norford, and W.~Kempton, ``{Shaving residential air-conditioner
		electricity peaks by intelligent use of the building thermal mass},''
	\emph{Energy}, vol.~16, no.~7, pp. 1001--1010, 1991.
	
	\bibitem{Roth2009}
	K.~Roth, J.~Dieckmann, and J.~Brodrick, ``{Emerging Technologies: Using
		Off-Peak Precooling},'' \emph{ASHRAE Journal}, vol.~51, no.~3, mar 2009.
	
	\bibitem{Tindemans2015}
	S.~H. Tindemans, V.~Trovato, and G.~Strbac, ``{Decentralized Control of
		Thermostatic Loads for Flexible Demand Response},'' \emph{IEEE Transactions
		on Control Systems Technology}, vol.~23, no.~5, pp. 1685--1700, 2015.
	
	\bibitem{Ihara1981}
	S.~Ihara and F.~Schweppe, ``{Physically Based Modeling of Cold Load Pickup},''
	\emph{IEEE Transactions on Power Apparatus and Systems}, vol. PAS-100, no.~9,
	pp. 4142--4150, 1981.
	
	\bibitem{Maurer1977}
	H.~Maurer, ``{On Optimal Control Problems with Bounded State Variables and
		Control Appearing Linearly},'' \emph{SIAM Journal on Control and
		Optimization}, vol.~15, no.~3, pp. 345--362, may 1977.
	
	\bibitem{Cesari1983}
	L.~Cesari, \emph{{Optimization—Theory and Applications}}.\hskip 1em plus
	0.5em minus 0.4em\relax New York, NY, USA: Springer New York, 1983.
	
	\bibitem{ASHRAE2016}
	ASHRAE, ``{Standard 90.4-2016 -- Energy Standard for Data Centers},'' 2016.
	
\end{thebibliography}

\appendix
\numberwithin{equation}{section}
\section{Description of the numerical method}\label{Section:Control-Parametrization}
We use the control parametrization method \cite{Teo1991,Rehbockt1999} to obtain an approximate solution to the optimal control problem. Let $\{t_{k}\}_{k = 0}^{n_{p}}$ denote a sequence of time points used to partition the control horizon $[0,T]$ into $n_{p} > 1$ subintervals, where $n_{p}$ is an integer, and let $\mathcal{U}_{n_{p}} \subset \mathcal{U}$ denote the corresponding subclass of step controls (cf. \eqref{eq:Piecewise-Constant-Control}). Each control $u \in \mathcal{U}_{n_{p}}$ is parametrized by an $n_{p}$-dimensional vector with components $u_{k} \in [0,1]$ for $k = 1,\ldots,n_{p}$.

The original optimal control problem is approximated by optimizing over the smaller subclass of controls $\mathcal{U}_{n_{p}}$, which can be treated as a constrained nonlinear optimization problem over the bounded $n_{p}$-dimensional parameter space defining controls $u \in \mathcal{U}_{n_{p}}$. Standard optimization packages can be used to solve this problem, and we used the Sequential Least Squares Programming (SLSQP) routine in Python.

\subsection{Loss functions for state and control constraints}\label{Section:Constraints}
In order to use the control parametrization method, we define constraints for the state and control variables in functional form. We only describe the constraints for Problem \eqref{Problem:Optimal-Control-Delivery-Problem}, noting that those for the other optimal control problems are defined analogously. First define $(t,x,u) \mapsto \psi_{1}(t,x,u)$ and $x \mapsto \phi_{1}(x)$ by,
\begin{equation}\label{eq:State-Constraints}
\begin{cases}
\psi_{1}(t,x,u) = (X_{max} - x)(x - X_{min}) \\
\phi_{1}(x) = (\hat{X} - x)(x - X_{min})
\end{cases}
\end{equation}
By definition, we say that the integral constraint $\psi_{1}$ is satisfied at $(t,x,u)$ if and only if $\psi_{1}(t,x,u) \ge 0$. Similarly, the terminal constraint $\phi_{1}$ is satisfied at $x$ if and only if $\phi_{1}(x) \ge 0$. Equation \eqref{eq:State-Constraints} corresponds to the time-dependent pure state constraints on the internal temperature over $[0,T]$ (cf. \eqref{eq:Bulding-Temperature-SLA}) and at time $T$ (cf. \eqref{eq:Bulding-Temperature-Terminal-Constraint}). In a similar way we define integral and terminal constraints for DTU,
\begin{equation}\label{eq:Control-Constraints}
\begin{cases}
\psi_{2}(t,x,u) = u - \sum_{i=1}^{N}\bigl[u_{ref}(t) + u_{i,ask}\bigr]\mathbf{1}_{[s_{i},e_{i})}(t) \\
\phi_{2}(x) = 0
\end{cases}
\end{equation}

Using these constraints we define {\it loss rate functions} $(t,x,u) \mapsto \Psi_{\eta}(t,x,u)$ and {\it terminal loss functions} $x \mapsto \Phi_{\eta}(x)$, $\eta \in \{1,2\}$, by,
\begin{equation}\label{eq:Lagrange-Constraints}
\begin{cases}
\Psi_{\eta}(t,x,u) = \left(\min(0,\psi_{\eta}(t,x,u))\right)^{2} \\
\Phi_{\eta}(x) = \left(\min(0,\phi_{\eta}(x))\right)^{2}
\end{cases}
\end{equation}
The loss functions $\Psi_{\eta}$ and $\Phi_{\eta}$ are combined to create a {\it total loss} for the constraints,
\begin{equation}\label{eq:Constraint-Canonical-Costs}
C_{\eta}(u) = \int_{0}^{T}\Psi_{\eta}(t,x(t),u(t))\,{d}t + \lambda_{\eta} \Phi_{\eta} (x(T)),
\end{equation}
where $\lambda_{\eta} > 0$ is a weighting for the terminal loss function $\Phi_{\eta}$. The total loss $C_{\eta}$ is non-negative by construction and is equal to zero if, equivalently, the relevant constraints are satisfied on $[0,T]$. We relax this condition by requiring,
\begin{equation}\label{eq:Canonical-Inequality-Constraints}
C_{\eta}(u) \le \varepsilon_{\eta},\; \eta \in \{1,2\},
\end{equation}
where $\varepsilon_{\eta} \ge 0$ is a sufficiently small tolerance parameter.

The SLSQP routine requires derivative information as input, but some of the constraints and costs in \eqref{Problem:Optimal-Control-Reference-Problem}, \eqref{Problem:Optimal-Control-Availability-Problem}, and \eqref{Problem:Optimal-Control-Delivery-Problem} are expressed in terms of indicator and ramp functions that are not smooth. Therefore, where necessary we approximate these functions smoothly as follows:
\begin{align*}
\mathbf{1}_{[0,\infty)}(y) & \approx \frac{e^{\theta y}}{1 + e^{\theta y}}, \quad \mathbf{1}_{[a,b]}(y) \approx \left(\frac{e^{\theta (y-a)}}{1 + e^{\theta (y-a)}}\right)\left(\frac{e^{\theta (b-y)}}{1 + e^{\theta (b-y)}}\right)\;\;\text{for}\;\; a < b, \\
\max(0,y) & \approx \frac{1}{\theta}\log(1 + e^{\theta y}),
\end{align*}
where $\theta > 0$ is a sufficiently large parameter.
\newpage
\section{Optimal power usage illustrations}\label{Appendix:Trajectories}

\begin{figure*}[!ht]
	\centering
	\subfigure[$\hat{X} = 18$ ($^{\circ}\rm{C}$)]{%
		\resizebox*{6cm}{!}{\includegraphics{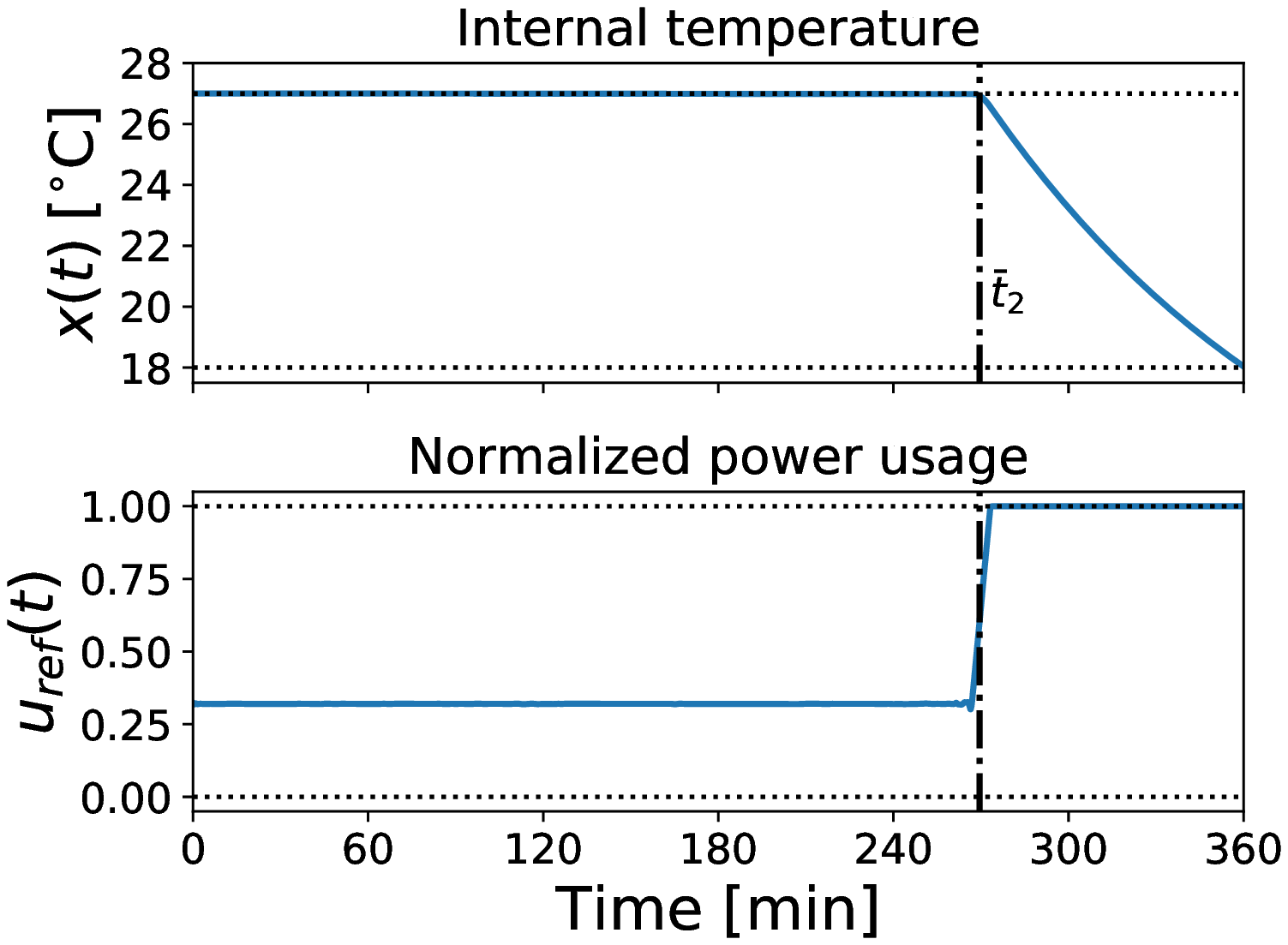}}
		\label{fig:Optimized-Reference-Power-1}
	}
	\hspace{5pt}
	\subfigure[$\hat{X} = 20$ ($^{\circ}\rm{C}$)]{%
		\resizebox*{6cm}{!}{\includegraphics{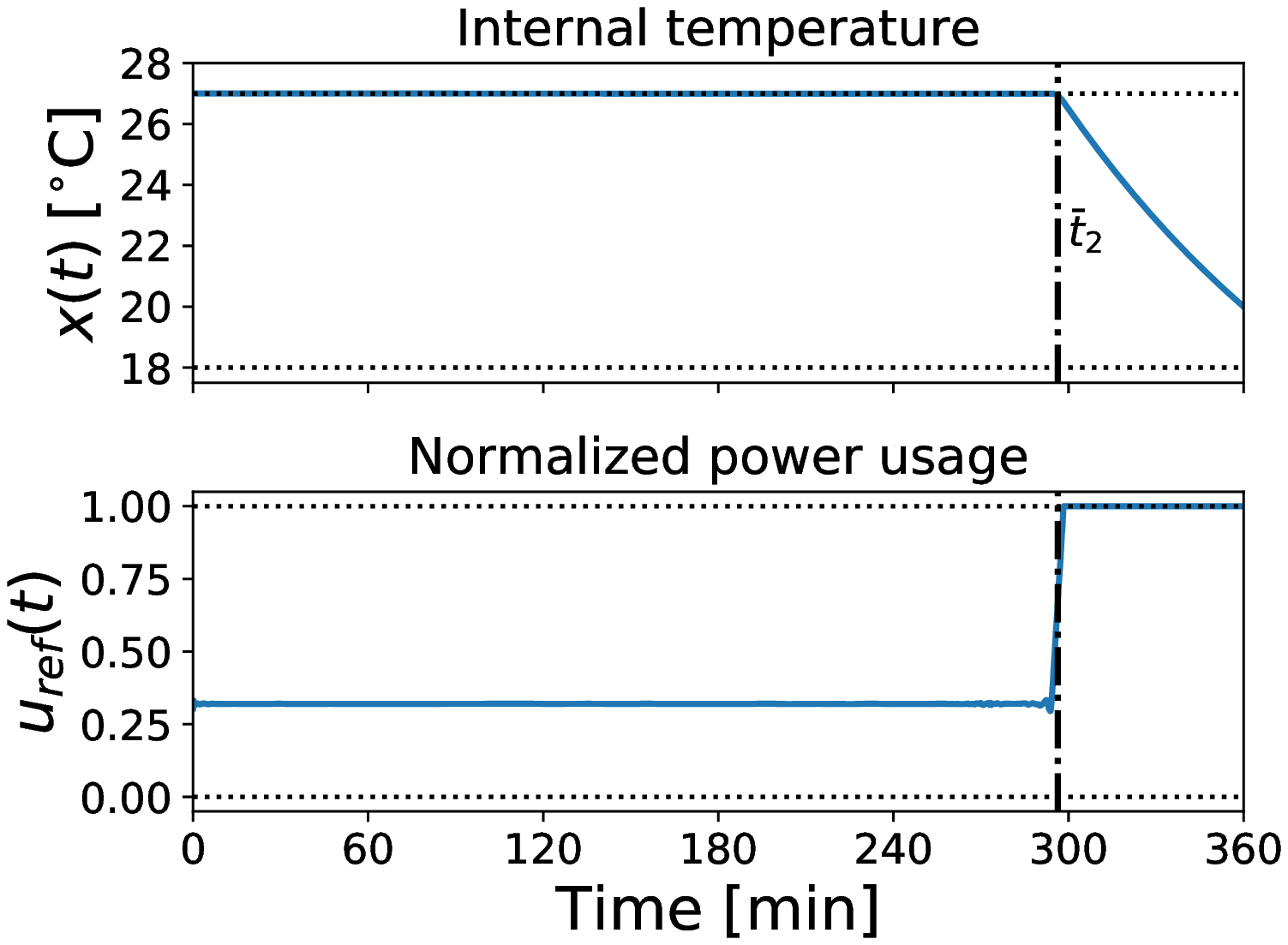}}
		\label{fig:Optimized-Reference-Power-2}
	}
	\caption{Optimized internal temperature $u$ and reference power profile $u_{ref}$ corresponding to two different values for the pre-cooling temperature $\hat{X}$. The profile $u_{ref}$ keeps the temperature constant at the maximum level for some time, and then quickly ramps up to full power around time $\bar{t}_{2}$ (vertical dash-dotted line, see \eqref{eq:Approximate-Full-Power-Time}) so that the temperature hits $\hat{X}$ exactly at the terminal time. The initial temperature is $x(0) = 27$ and regularizer has value $\alpha_{ref} = 0.01$. Control and temperature constraints are shown using dotted horizontal lines.}
	\label{fig:Optimized-Reference-Power}
\end{figure*}

\begin{figure*}[!ht]
	\centering
	\subfigure[$\hat{X} = 18$]{%
		\resizebox*{6cm}{!}{\includegraphics{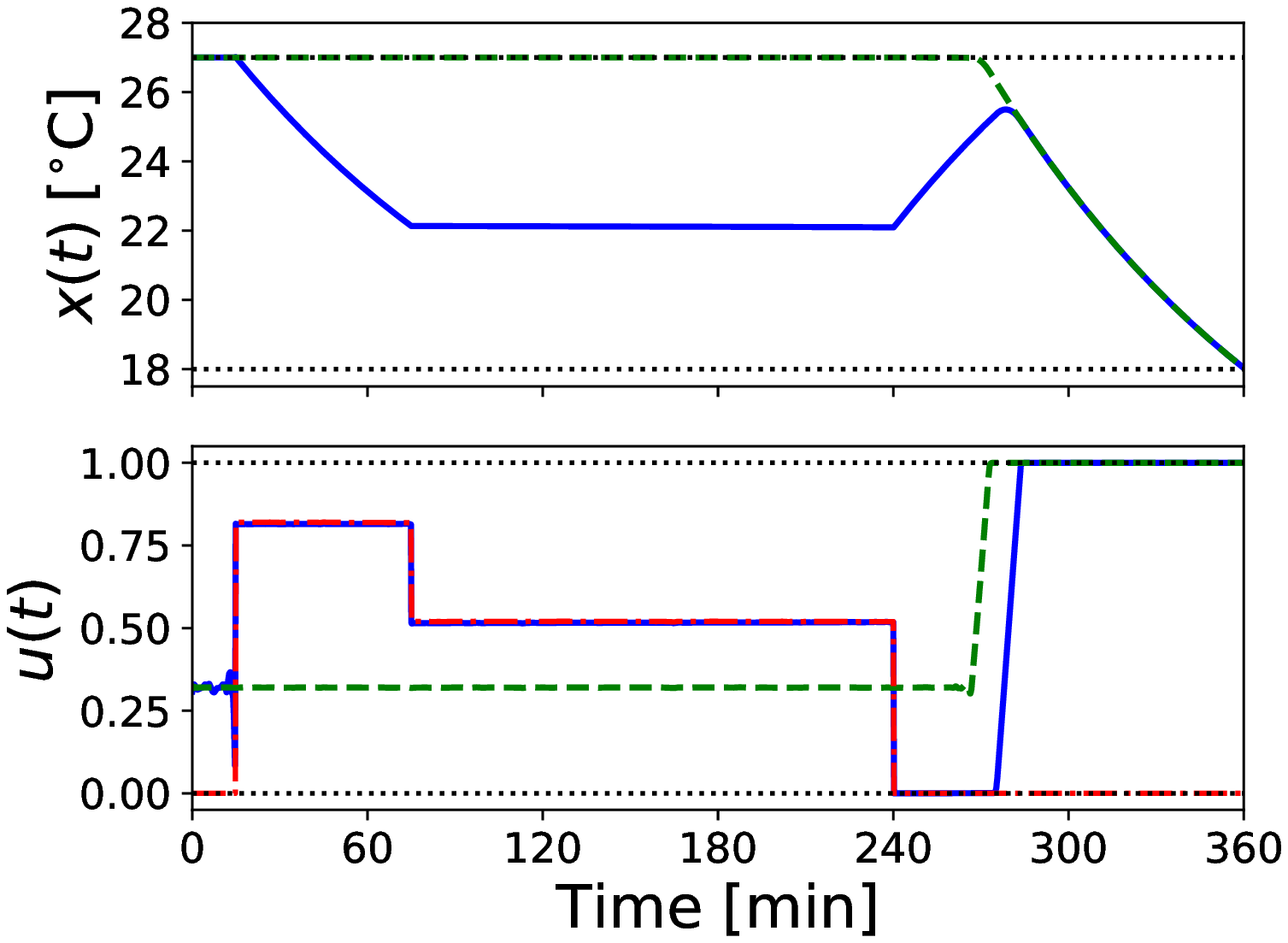}}
		\label{fig:Optimized-Delivery-Power-1}
	}
	\hspace{5pt}
	\subfigure[$\hat{X} = 20$]{%
		\resizebox*{6.15cm}{!}{\includegraphics[height=0.22\textheight, width=0.475\textwidth]{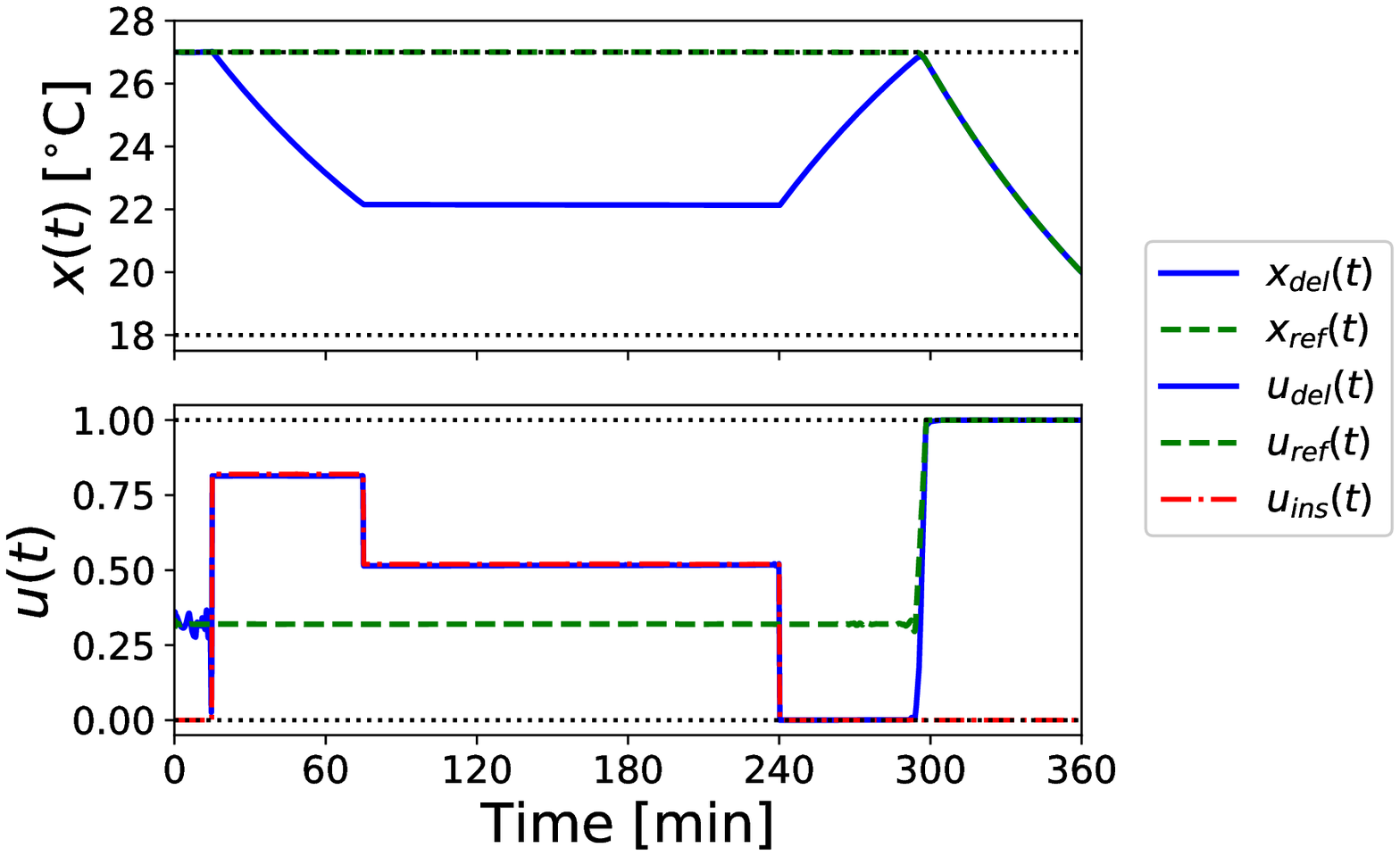}}
		\label{fig:Optimized-Delivery-Power-2}
	}
	\caption{Optimized internal temperature $x$ and delivery, reference, and minimum power profiles, $u_{del}$, $u_{ref}$, and $u_{ins}$ corresponding to two different values for the pre-cooling temperature $\hat{X}$, shown using the dot-dash line. Reserve service instructions call for an increase in normalized power by $0.5$ between minutes 15 and 75, then $0.2$ between minutes 75 and 240. The profile $u_{del}$ uses minimal power to satisfy the temperature constraints and reserve service instructions. The building is also pre-cooled after the service has been delivered. Parameters values are $x(0) = 27$ for the initial temperature and $\alpha_{ref} = \alpha_{del} = 0.01$ for the regularizers. Temperature and control constraints are shown using dotted horizontal lines.}
	\label{fig:Optimized-Delivery-Power}
\end{figure*}

\begin{figure*}[!ht]
	\centering
	\subfigure[$\alpha_{alt} = 10$]{%
		\resizebox*{6cm}{!}{\includegraphics{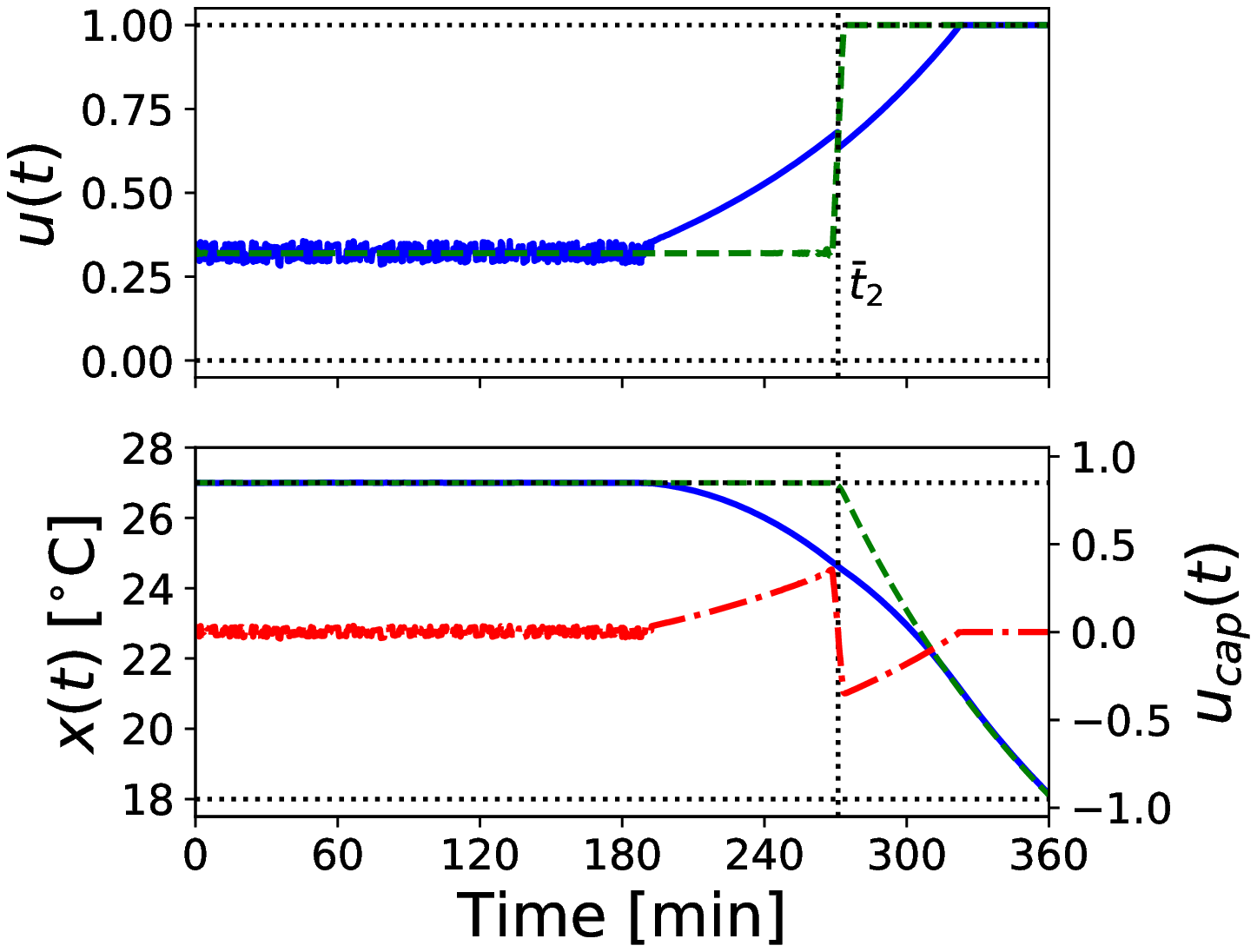}}
		\label{fig:Optimized-Alternative-Power-4}
	}
	\hspace{5pt}
	\subfigure[$\alpha_{alt} = 1$]{%
		\resizebox*{6cm}{!}{\includegraphics{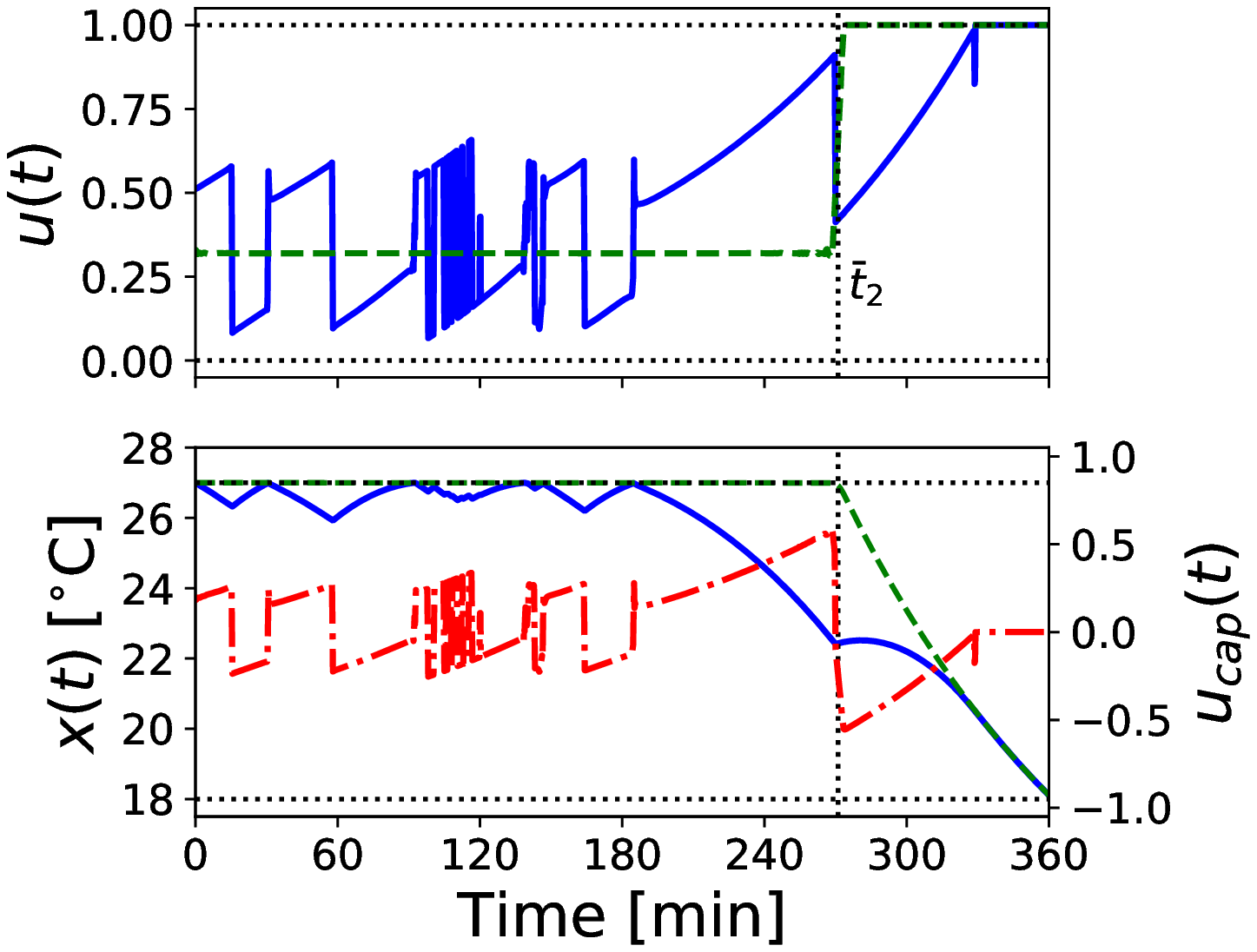}}
		\label{fig:Optimized-Alternative-Power-5}
	}
	\hspace{5pt}
	\subfigure[$\alpha_{alt} = \frac{1}{10}$]{%
		\resizebox*{6.15cm}{!}{\includegraphics{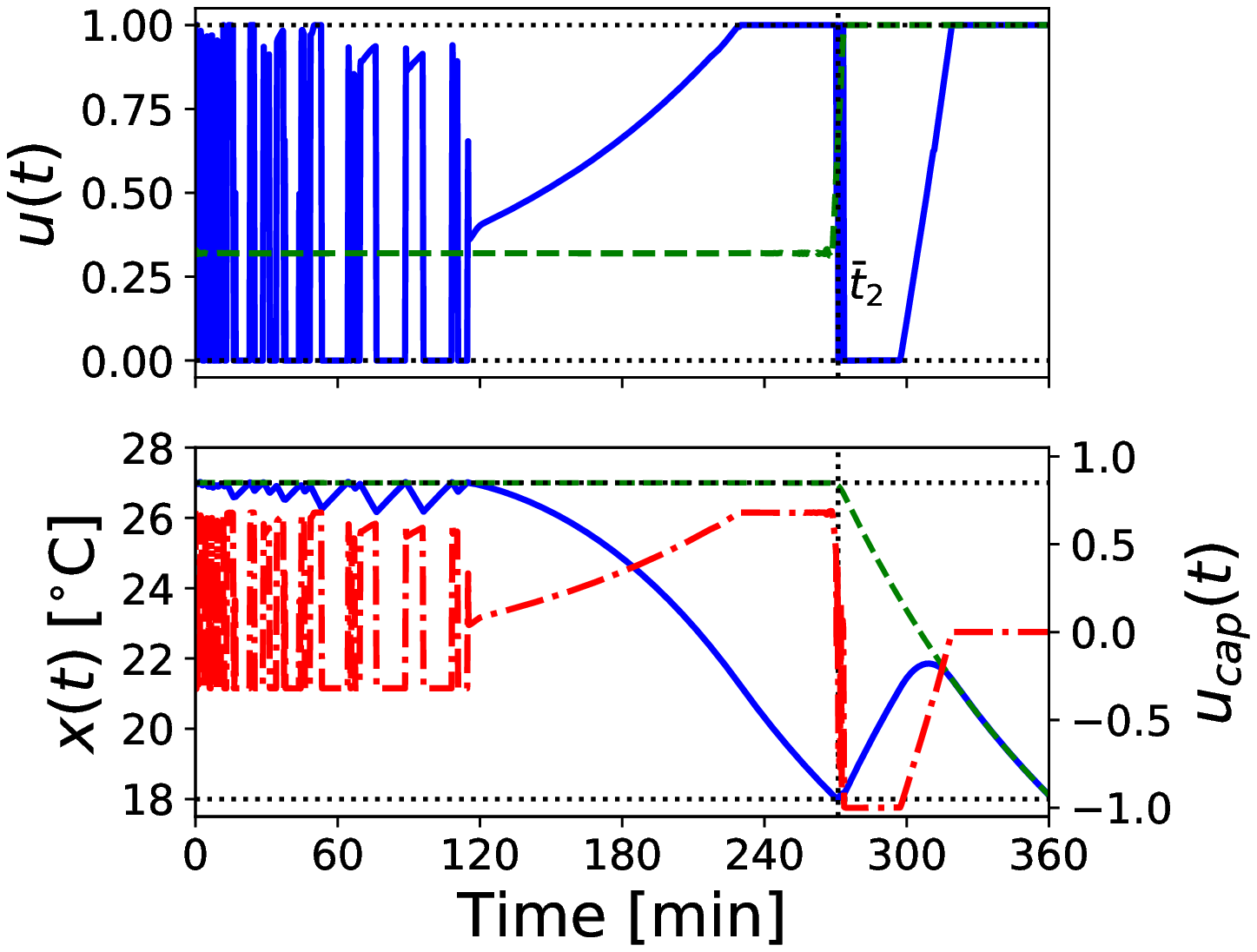}}
		\label{fig:Optimized-Alternative-Power-6}
	}
	\hspace{5pt}
	\subfigure[$\alpha_{alt} = \frac{1}{100}$]{%
		\resizebox*{6.35cm}{!}{\includegraphics[height=0.225\textheight, width=0.475\textwidth]{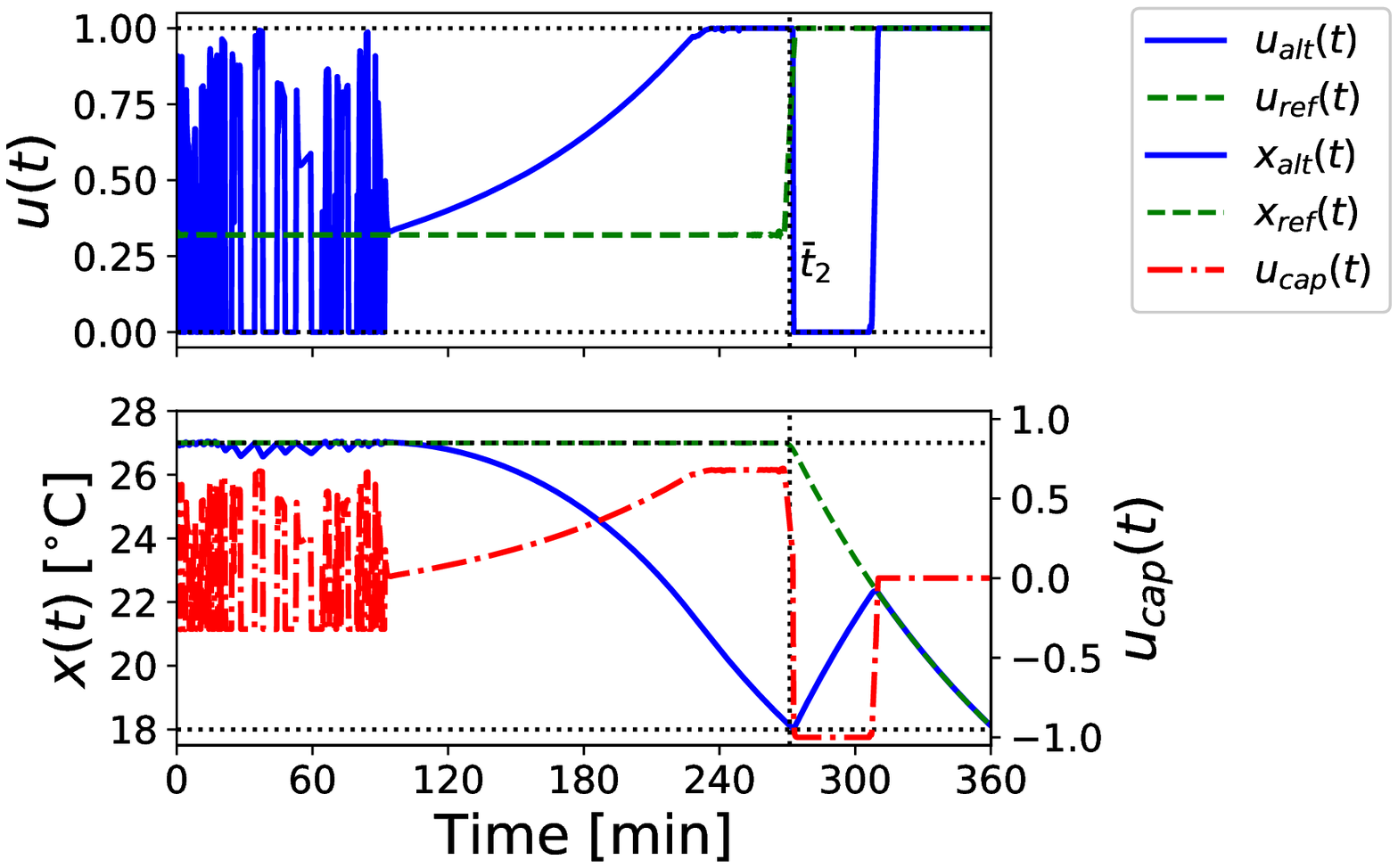}}
		\label{fig:Optimized-Alternative-Power-7}
	}
	\caption{Optimized internal temperature $x$, alternative power usage $u_{alt}$ and reserve capacity $u_{cap}$ corresponding to the benefit-cost ratio $\frac{R}{P} = 1$ and decreasing values for the regularizer $\alpha_{alt} \in \{10,1,\tfrac{1}{10},\frac{1}{100}\}$. When $\alpha_{alt}$ is large, the alternative profile $u_{alt}$ tries to minimize the overall cost of consumption, similarly to the reference profile $u_{ref}$. However, as $\alpha_{alt}$ decreases $u_{alt}$ increasingly tries to stay below $u_{ref}$, leading to rapid bang-bang control behaviour. Parameter values are $x(0) = 25$ and $\hat{X} = 18$ for the initial and pre-cooling temperatures respectively. Temperature and control constraints are shown using dotted horizontal lines whilst $\bar{t}_{2}$ (cf. \eqref{eq:Approximate-Full-Power-Time}) is shown using the dash-dotted vertical line.}
	\label{fig:Optimized-Alternative-Power-Equal-Reward-Cost}
\end{figure*}

\end{document}